%% file: version-1.0.tex
\documentclass[12]{siamltex}
\usepackage{amsmath,amssymb,amsfonts,graphicx,fancybox,shadow,color}
\usepackage{hyperref}
\usepackage{wrapfig}
\newtheorem{remark}[theorem]{Remark}

\def\be#1\ee{\begin{equation}#1\end{equation}}

\input{math-commands.tex}

\newcommand{\RR}{{\mathbb{R}}}

\newcommand{\tilV}{{\tilde{V}}}

\def\tF {\tilde{F}}
\def\tM {\tilde{M}}
\def\tS {\tilde{S}}
\def\tc {\tilde{c}}
\def\tC {\tilde{C}}
\def\tb {\tilde{b}}
\def\tB {\tilde{B}}

\def\uu{{\mathbf u}}

\def\be#1\ee{\begin{equation}#1\end{equation}}

\def\bF{{\mathbf F}}

\begin{document}
\pagestyle{myheadings}

\title{ Regularized Reduced Order Lippmann-Schwinger-Lanczos Method for Inverse Scattering Problems in the Frequency Domain 
}
\author{
 J. Baker\footnotemark[1],  E. Cherkaev\footnotemark[2], V. Druskin\footnotemark[3], S. Moskow\footnotemark[4] and M. Zaslavsky\footnotemark[5]}

\renewcommand{\thefootnote}{\fnsymbol{footnote}}

\footnotetext[1]{Department of Mathematics, University of Utah, 155 South 1400 East, JWB 233, Salt Lake City, UT 84112 (baker@math.utah.edu)}
\footnotetext[2]{Department of Mathematics, University of Utah, 155 South 1400 East, JWB 233, Salt Lake City, UT 84112 (elena@math.utah.edu)}

\footnotetext[3]{Department of Mathematical Sciences, Worcester Polytechnic Institute, Stratton Hall,
100 Institute Road, Worcester MA, 01609 (vdruskin1@gmail.com)}
\footnotetext[4]{Department of Mathematics, Drexel University, Korman Center, 3141 Chestnut Street, Philadelphia, PA 19104
(moskow@math.drexel.edu)}
\footnotetext[5]{Department of Mathematics, Southern Methodist University, 6425 Boaz Lane, Dallas TX 75205  (mzaslavsky@slb.com)}

\maketitle
\begin{abstract}  
Inverse scattering has a broad applicability in quantum mechanics, remote sensing, geophysical, and medical imaging. This paper presents a robust direct reduced order model (ROM) method for solving inverse scattering problems based on an efficient approximation of the resolvent operator regularizing the Lippmann-Schwinger-Lanczos (LSL) algorithm. 
We show that the efficiency of the method relies upon the weak dependence of the orthogonalized basis on the unknown potential in the Schrödinger equation by demonstrating that the Lanczos orthogonalization is equivalent to performing Gram-Schmidt on the ROM time snapshots. 
We then develop the LSL algorithm in the frequency domain with two levels of regularization. 
We show that the same procedure can be extended beyond the Schrödinger formulation to the Helmholtz equation, e.g., to imaging the conductivity using diffusive electromagnetic fields in conductive media with localized positive conductivity perturbations. 
Numerical experiments for Helmholtz and Schrödinger problems show that the proposed bi-level regularization scheme significantly improves the performance of the LSL algorithm, allowing for good reconstructions with noisy data and large data sets.
\end{abstract}

%


\section{Introduction}

The inverse scattering problem formulated for the Helmholtz and Schrödinger operators arises in various fields, including quantum mechanics, remote sensing, geophysical, and medical imaging.
The goal of imaging is to find the potential or properties of the medium in the domain using near-field measured data. There are different approaches to the inverse problem, starting with an original method based on the solution of a set of integral equations developed in celebrated works of Gelfand \& Levitan, Marchenko, and Krein. 
 In different formulations appropriate for particular applications, the existence and uniqueness of the solution of inverse scattering problems, as well as computational approaches, were discussed in many now classical works.  
 Efficient numerical methods for inverse scattering have been developed in application to acoustic imaging, electromagnetic sensing, and seismic exploration; among them are iterative methods for full-wave inversion based on the adjoint or backpropagation method, techniques based on Born and Rytov approximations, layer stripping methods, asymptotic methods specifically developed for small volume inhomogeneities, as well as Kirchhoff migration and solving the Lippmann-Schwinger (LS) integral equation, see for instance  
\cite{nachman1984multidimensional,beylkin1990linearized,chew1990reconstruction,habashy1993beyond,Biondi,vogelius2000asymptotic,bleistein2001mathematics,ammari2003boundary,cheney2009fundamentals, CakHad2,symes2007reverse,KiMoSc,colton2013inverse,ammari2012direct}
  and references therein. 
 Formulated as the identification of the parameters of a layered medium, inverse scattering is related to cascade realizations of the transfer functions for digital filters or the minimal realization problem, actively discussed in electrical engineering and control theory literature \cite{gragg1983partial,bruckstein1988inverse,parlett1992reduction}. 
 Following the terminology used in that field, we refer to the inverse problems as single-input single-output (SISO) and multiple-input multiple-output (MIMO), depending on the problem setting.

The present work develops a robust direct reduced-order model method for solving inverse scattering problems for the Schrödinger and Helmholtz equations based on a Lippmann-Schwinger-Lanczos (LSL) algorithm in the frequency domain with two levels of regularization. This is the third paper in a series on the reduced order model (ROM)  Lippmann-Schwinger-Lanczos (LSL) algorithm 
\cite{druskin2021lippmann},\cite{druskin2022extension}. 
Here, we advance the algorithm foundation,  develop a robust regularization scheme, and extend the area of applicability of the method.  
Presented numerical experiments demonstrate that the regularized LSL method results in accurate reconstructions with noisy data and large data sets.

Our approach is based on the approximation of a resolvent operator and constructing a Reduced Order Model (ROM) that allows effectively substitute solving a large-scale inverse problem with a problem of much smaller size. Such a reduction of the size of the numerical problem has irrefutable advantages in saving computational time and required storage memory. ROMs are generally constructed by projecting the solution of a large-scale problem into an appropriate subspace \cite{Antoulas01asurvey,
willcox2002balanced,guttel2013rational,benner2015survey}; this can be achieved by different methods such as Krylov subspaces, proper orthogonal decomposition, truncated SVD, etc.  
 Here, we use a truncated Lanczos representation. Lanczos method is intrinsically related to Krylov subspace methods efficiently used in the iterative algorithms \cite{feldmann1995efficient,grimme1997krylov,
zimmerling2016lanczos,kordy2016adaptive,kordy2017null}; the truncation results in a reduced order solution spanned by the dominant eigenfunctions of the operator \cite{benner2014model,zimmerling2017model,cherkaev1996inverse,baker2023learning}.  
There is extensive literature discussing these approaches for particular applications; see \cite{BAI2002,antoulas2005approximation,benner2014model,guttel2013rational} and references therein.


Our contribution is three-fold. First, we investigate the weak dependence of the Lanczos orthogonalized frequency domain snapshots on perturbations of the potential
by viewing the ROM as a spectral projection of the time-domain problem; the well-known interpolatory-projection dualism \cite{BEATTIE20155}.
This allows us to invoke the reflection cancellation properties of the time-domain snapshots \cite{druskin2016direct} to explain the weak dependence in the frequency domain.

Next, we consider LSL regularization. Generally, the multidimensional inverse scattering problem is ill-posed and requires regularization to develop a stable numerical algorithm. This ill-posedness is aggravated if the scattering data is obtained in the diffusive regime, as in, for example, diffuse optical tomography or control source electromagnetic exploration. 
In the LSL method, this ill-posedness creates two sources of instability. One is the conventional instability of the Lippmann-Schwinger (LS) integral equation, which, as a linear Fredholm integral equation, can be resolved by SVD truncation or other regularization techniques developed for linear problems \cite{vogel2002computational,bertete2002non,kilmer2007projection,hansen2010discrete,reichel2012tikhonov,chung2015framework}. We use a truncation of the spectral decomposition in the current work to deal with this type of instability.
Among various regularization methods suggested for the solution of linear and non-linear ill-posed inverse problems, truncation of a spectral decomposition allows for a significant reduction in the size of the computational problem, which is paramount for large-scale problems. 
 
The other stability problem is the ill-conditioning of the projection subspace. 
The developed method constructs the data-driven ROM and Gramian directly from the data (that is why they are called data-driven); the ill-conditioning   
 may lead to the ROM's sensitivity to data errors and may cause nonphysical indefiniteness of the mass and stiffness matrices.  In turn, this leads to the loss of the Hamiltonian property of the
underlying reduced order dynamical system. To circumvent this problem, a ROM regularization via data-driven Gramian truncation was introduced in \cite{borcea2019robust, borcea2020reduced} for inverse scattering problem in the time domain. In those works, the data-driven ROM constructed in the time domain was applied respectively to the so-called data-to-born transform (removing multiple echoes directly from the data) and to preconditioned full wave inversion. In both cases, regularization was performed via  SVD truncation of the data-driven approximate Gramian to the level of the measurement error. In the current paper, we extend this approach to the frequency domain LSL formulation. An advantage of LSL is that it allows for the two levels of truncation. 
Thanks to this, we can relax the truncation level of the Gramian, benefiting the inversion quality. We should also mention that the approximate data-driven Gramian truncation developed here can be viewed as a particular case of the recently emerged data-driven balanced truncation method \cite{gosea2022data}, although, for LSL, its implementation requires some specifics.

Finally, we extend the LSL approach from the Schrödinger problem to the \\Helmholtz variable coefficient formulation. Such an extension allows us to broaden the application of LSL, for example, to low-frequency quasi-stationary electromagnetic problems with variable conductivity. We can qualitatively (and in some cases quantitatively) map conductivity perturbations, provided they do not significantly affect the high-frequency asymptotics of the host models, as in the case of localized positive perturbations. Extensions to a more general class of Helmholtz problems could be possible with additional iterations; see \cite{doi:10.1137/22M1517342} for the Lippmann-Schwinger formulation in the time-domain, however adding iterations would eliminate the great advantage of the direct inversion method.  

The paper is organized as follows. We describe the LSL method for a SISO problem in Section \ref{reviewSISO} and summarize the steps of the LSL algorithm in Section \ref{LSLsteps}.  
In Section \ref{ssec:time}, we discuss the connection between the frequency domain LSL and time domain snapshots of the solution; in particular, we show that Lanczos orthogonalization corresponds to Gram-Schmidt on the corresponding sequential time snapshots. 
The regularized Lippmann-Schwinger-Lanczos method is presented in Section \ref{reg-LSL}. 
Extensions to the multidimensional MIMO problem and to the Helmholtz equation modeling electromagnetic diffusion are described in 
Section \ref{MIMO-LSL} and in Section \ref{H-LSL}, respectively.  Results of numerical experiments are shown in Section \ref{Num}.



\section{The LSL method for SISO problems} \label{reviewSISO}

We give a review of the method and some necessary results from
\cite{druskin2021lippmann}. 
We start with considering the single input single output (SISO) inverse problem in one dimension. The forward problem in the frequency domain is described by the Schrödinger equation with a potential $p(x)$ for a scalar function $u$ 
	\be\label{eq:1D} 
	-u^{\prime\prime}(x) +p(x)u(x)+\lambda u(x)=g(x) \quad \mbox{in} \,\,  \Omega ,
	 \quad u^\prime(0)=0,\ u^\prime(L) =0,
	\ee
where  $\Omega =(0,L)$ and $0<L\le\infty$. The source $g(x)$ is assumed to be a compactly supported real distribution localized near the origin, for example, $ g = \delta(x-\epsilon)$ with small $\epsilon>0$, or a slightly smoothed approximation. Note that when $g$ is a delta function at the origin, this corresponds to an inhomogeneous Neumann boundary condition. We can write the solution formally using the resolvent
	$$
	\label{eq:resolvent} 
	u= \left({\cal L}+\lambda I \right)^{-1} g, \qquad  
\mbox{ where} \qquad  
{\cal L}=- \Delta  
+{p}I, 
$$ 
and where the inverse is understood to correspond to imposing homogeneous Neumann boundary conditions.

 We note that $u$ depends on $\lambda$, but we will denote it explicitly only when necessary for clarity. 
 It is known that for nonnegative $p$, the above resolvent is well defined for $\lambda$ that is not on the negative real axis.  The SISO transfer function is defined as 
\be\label{eq:transfer}
	F(\lambda)=\int_0^L g(x) u(x;\lambda)dx =\langle g,u\rangle =
	\langle g,\left({\cal L}+\lambda I \right)^{-1}g\rangle 
	\ee
The notation  $\langle, \rangle$ is used throughout the paper to denote the continuous $L^2(\Omega)$ inner product $$\langle w,v\rangle=\int_\Omega {w}(x)v(x)dx ,$$ 
where  $\Omega =(0,L)$ in the one dimensional case. 
The data are the values of the function $F(\lambda)$ and its derivative given for specific $\lambda_j\in \RR $, $j=1,\ldots, m$: 
\begin{equation} \label{realdata} 
F(\lambda)|_{\lambda=\lambda_j}\in\RR, \ \ \ \frac{dF(\lambda)}{d\lambda} |_{\lambda=\lambda_j}\in \RR \ \ \ \mbox{for} \ \ j=1,\ldots, m.
\end{equation} 
The SISO inverse problem requires to determine $p(x)$ in (\ref{eq:1D}) from the data (\ref{realdata}). 
Such a formulation can be obtained, for example, using Zolotarev-optimal \cite{druskin2009solution} and data-driven $H_2$-optimal ROMs \cite{BEATTIE20155} for a diffusion problem with $\frac{\partial} {\partial t}$ in place of $\lambda$ in (\ref{eq:1D}). It is known that such approximations give the most accurate solutions to the inverse problem \cite{druskin2013solution}, although their data-driven construction presents significant computational challenges, especially for MIMO problems. Moreover, their stability for noisy data is not well-studied. 
We also note that the following discussion can be extended for complex data points
outside of the spectrum of $\cal{L}$ (see \cite{druskin2021lippmann} for a more thorough treatment of the case of complex data points).

The Lippmann-Schwinger (LS) integral equation is fundamental to many inverse scattering algorithms.
If  $F_p$ denotes the data (\ref{eq:transfer}) for the PDE with the unknown potential  $p$, and 
 $F_0$ is the background transfer function corresponding to the PDE with a known $p_0$ (here we assume $p_0=0$), then 
the Lippmann-Schwinger equation for the unknown coefficient  $p$ can be written as 
\begin{equation}\label{eq:LipSwi} 
F_p -F_0 =-\langle u^0,p u_p\rangle 
\end{equation}
where 
$u_p$ is  the   internal PDE solution   corresponding to the unknown coefficient $p$, and $u^0$ is the background internal PDE solution. It is assumed that $u^0$ can be computed by solving the  PDE with known background parameters, while  $u_p$ is unknown because it depends on the unknown $p$. This dependence makes  \eqref{eq:LipSwi} nonlinear, which is one of the main difficulties of applying the  LS approach to inverse problems. A conventional method \cite{chew1990reconstruction} is the so-called distorted Born iteration, i.e., to use the background solution in place of $u_p$ on the first iteration and successively update by using the approximation to $p$ obtained in the previous iteration as the background medium. This procedure can be quite computationally expensive and is very sensitive to the initial guess. 

To circumvent these problems, 
the Lippmann-Schwinger-Lanczos method was introduced in \cite{druskin2021lippmann,druskin2022extension}. 
The LSL method uses the data-driven internal solution ${\bf u}_p$ in place of $u_p$:
\begin{equation}\label{eq:LipSwiL} 
F_p -F_0 \approx -\langle u^0,p {\bf u}_p\rangle.
\end{equation}
Solution ${\bf u}_p$ can be computed directly from the data without knowing $p$; then 
(\ref{eq:LipSwiL}) becomes linear with respect to $p$. This precomputing is based on embedding properties of the data-driven reduced order models (ROMs) developed in \cite{druskin2016direct}. The LSL approximate of the internal solution ${\bf u}_p$ is computed via Lanczos orthogonalization, hence the name.


Let $u_j = u(x,\lambda_j)$ be solutions to (\ref{eq:1D}) corresponding to $\lambda=\lambda_j$ for $j=1,\ldots, m$. We consider the projection given by the matrix \footnote{ Basis $V$ can be viewed as a matrix with infinite columns or equivalently as a row-vector of continuous functions, in either case, we use the $L_2$ inner product  $\langle,\rangle$ defined earlier.} 
$$
V=\left[ u_1(x),\ldots , u_m(x) \right]\in\RR^{\infty\times m}.
$$
The Galerkin system (\cite{fletcher1984computational}) projecting the problem (\ref{eq:1D}) into the subspace $V$ spanned by the functions  $ u_1(x),\ldots , u_m(x)$, 
\be \label{eq:ROM}
(S +\lambda M)c =b, 
\ee
determines the data-driven Reduced Order Model (ROM). 
Here $S, M \in \RR^{m\times m}$  in (\ref{eq:ROM}) are symmetric, positive definite stiffness and mass matrices, respectively, given by 
$$
S_{ij}=\langle u_i, {\cal L} u_j \rangle , 
\qquad
M_{ij}=\langle {u}_i, u_j \rangle. 
$$ 
The right-hand side $b\in\RR^m$ in (\ref{eq:ROM}) is a column vector with components 
$$b_i =\langle  {u}_i, g  \rangle. $$ 
The Galerkin solution for the system (\ref{eq:ROM})  is determined by the vector-valued function of $\lambda$,  $c \in\RR^{m}$, where $c$  corresponds to a column vector of coefficients with respect to the basis $V$ of exact solutions. For any $\lambda$, the solution to (\ref{eq:1D}) can be approximated by its Galerkin projection 
\[
u\approx \hat{u} ={V}c ={V} (S+\lambda M)^{-1}b. 
\] 
A key component of the data-driven ROM approach is that even though the solutions $u$ are unknown, the matrices $S$ and $M$ can be obtained from the data.

\section{The data driven LSL algorithm} \label{LSLsteps}

In the data driven algorithm, the mass and stiffness matrices $M$ and $S$ are computed directly (\cite{fletcher1984computational,BEATTIE20155,BoDrMaMoZa}) from the data as: 
\begin{equation}
\label{eq:massmtr}
M_{ij}=\frac{{F}(\lambda_i)-F(\lambda_j)}{\lambda_j-{\lambda}_i}, \ \ \ M_{ii} = -\frac{dF}{d\lambda}(\lambda_i).
\end{equation}
and
\begin{equation}
\label{eq:stifmtr}
S_{ij}=\frac{F(\lambda_j)\lambda_j-{F}(\lambda_i){\lambda}_i}{\lambda_j-{\lambda}_i}, \ \ \ S_{ii} = \frac{d(\lambda F)}{d\lambda}(\lambda_i).
\end{equation}
Indeed, multiplying by $u_i$ equation (\ref{eq:1D}) written for $\lambda = \lambda_j$ and integrating by parts, and then switching $i$ and $j$, we obtain
\be \label{eq:SMmtr}
(S_{ij}+\lambda_j M_{ij}) = F(\lambda_i), \qquad
(S_{ij}+\lambda_i M_{ij}) = F(\lambda_j).
\ee
Subtracting the second equation from the first one and dividing by $(\lambda_j-{\lambda}_i)$ we obtain 
the elements of the matrix $M$  in (\ref{eq:massmtr}). Similarly, one obtains formulas for the elements of the stiffness matrix (\ref{eq:stifmtr}) by multiplying equations in (\ref{eq:SMmtr}) by $\lambda_i$ and $\lambda_j$ before doing the subtraction step. 
The LSL algorithm corresponds to executing the steps below.

\begin{enumerate}
   \item  {\bf Data generated ROM.}
The first step of the data driven algorithm is to compute  
the mass and stiffness matrices $M$ and $S$ (\cite{fletcher1984computational,BEATTIE20155,BoDrMaMoZa}) from the data using formulas (\ref{eq:massmtr}), (\ref{eq:stifmtr}) and the right-hand side vector $b\in\RR^m$  in (\ref{eq:ROM}). 
The ROM transfer function corresponding to the Galerkin system is   
     \[ \hat{F}(\lambda) := \langle \hat{u} , g\rangle =b^\top c . \] 
 It can be shown that due to the data matching conditions, this projected into the Galerkin subspace ROM transfer function matches the data exactly. 

\begin{proposition} \cite{BoDrMaMoZa,druskin2021lippmann}  \label{prop:1}
Assume that $M$ and $S$ are generated from the data  (\ref{realdata})  by the formulas (\ref{eq:massmtr}) and (\ref{eq:stifmtr}).  Then the Galerkin projection of the solution of (\ref{eq:1D})
  \[\hat{u}(\lambda)=Vc(\lambda)=V(S+\lambda M)^{-1}b\] 
is exact at $\lambda=\lambda_j$ , 
$$ \hat{u}(\lambda_j) = u(\lambda_j),$$
and hence
$$ \hat{F}(\lambda_j)=b^\top(S+\lambda_j M)^{-1}b=F(\lambda_j).$$
Furthermore
$$\frac{d\hat{F}}{d\lambda} (\lambda_j) =\frac{d F}{d\lambda} (\lambda_j)  $$  
for $j=1,\ldots, m$.
\end{proposition}
See \cite{BoDrMaMoZa} for the proof for real $\lambda$, and \cite{druskin2021lippmann}
for the corresponding result and proof for $\lambda_j$ complex.  

\item {\bf Lanczos orthogonalization.} The next step is to change the basis by orthogonalizing matrix $A=M^{-1}S$ using the Lanczos algorithm. More precisely, we run $m$ steps of the { $M$-symmetric}  Lanczos algorithm applied to  matrix $A=M^{-1}S$ and initial vector  $M^{-1}b$. This yields tridiagonal matrix $T\in\RR^{m\times m}$ and $M$-orthonormal Lanczos vectors $q_i\in \RR^m$, such that
	\be\label{eq:lancz} AQ =Q T , \qquad Q^\top MQ=I,\ee
where $$Q=[q_1, q_2, \ldots, q_m]\in{\RR^{m\times m}},$$ and $$q_1=M^{-1}b/\sqrt{b^\top M^{-1}b}.$$
 The new basis is orthonormal in $L^2(0,L)$ and is given by the row vector of continuous functions  $${VQ=  [ \sum_{j=1}^mq_{j1}u_j ,\ldots , \sum_{j=1}^m q_{jm}u_j ]} \in \mathbb{R}^{\infty\times m} .$$   The Galerkin solutions $\hat{u}$ and transfer function $\hat{F}(\lambda)$ can be written in this new basis as
\be \label{eq:state} 
\hat{u}(\lambda)=\sqrt{b^\top M^{-1}b}V Q(T+\lambda I)^{-1}e_1, 
\ee 
\be \hat{F}(\lambda)= (b^\top M^{-1}b)  e_1^\top (T+\lambda I)^{-1}e_1
\ee
where $e_1 = (1,0,\ldots,0)^T $ is the first coordinate column vector in $\RR^m$.

\item {\bf Internal solutions.} The orthogonalized ROM is used to produce internal solutions. 
As the true potential $p$ and the basis  $V$ consisting of exact solutions are not known, 
we replace the unknown orthogonalized internal solutions $VQ$ with orthogonalized background solutions $V_0Q_0$ corresponding to background $p_0=0$. Here $V_0$ is the row vector of background solutions 
$$V_0 = [ u^0_1,\ldots , u^0_m ] $$ 
to (\ref{eq:1D}) corresponding to $p=p_0=0$ and the same spectral points 
$\lambda= \lambda_1, \ldots \lambda_m$. A ROM for the background problem is computed in the same way, and $Q_0$ is computed using Lanczos orthogonalization.  The approximation  \be\label{eq:basis} VQ\approx V_0Q_0\ee is the crucial step -  
it is used in (\ref{eq:state}) to compute an approximation $\uu$ to the unknown internal solution 
$u(x,\lambda)$ 
\begin{equation}\label{internal1d} 
 u \approx \uu= \sqrt{b^\top M^{-1}b} V_0 Q_0(T+\lambda I)^{-1}e_1 .  
\end{equation}

\item {\bf Lippmann-Schwinger.} In the Lippmann-Schwinger formulation, one needs to solve the nonlinear inverse problem. However, the data driven approach allows us to reformulate it as a linear problem. 
From  (\ref{eq:transfer}) and its background counterpart, we obtain the Lippmann-Schwinger equation 
	\be  \label{eq:int1d}F_0(\lambda_j)-F(\lambda_j)=\int u^{0}_{j}(x) u_j(x)p(x)dx \ee
for $j=1,\ldots, m$. 
Here $F_0$ is the  transfer function corresponding to the background problem
\begin{equation} \label{eq:backtransfer}
F_0(\lambda) =\int_0^L g(x) u^0(x;\lambda )dx . \end{equation} 
For real $\lambda_j\in \RR$ we have the following $2m$ equations 
\begin{equation}\label{eq:intd}
F_0(\lambda_j) -F(\lambda_j) =\int  u^0_j(x) u_j(x)p(x)dx 
\end{equation} 
and 
\begin{equation}  
\frac{d}{d\lambda}(F_0-F)|_{\lambda= \lambda_j}=\int \frac{d}{d\lambda} [u^0(x;\lambda) u(x;\lambda)]_{\lambda=\lambda_j} p(x)dx,  
\label{eq:intdd} 
\end{equation}
for $j=1,\ldots,m$. 
The internal solutions $u_j(x)$ and their derivatives with respect to $\lambda$ are unknown and depend on unknown $p$, which makes the LS problem (\ref{eq:int1d}) nonlinear. To reduce it to a linear problem  we replace $u$  in (\ref{eq:intd}),(\ref{eq:intdd}) with its approximation $ u \approx \uu $   computed using  
 (\ref{internal1d}) \cite{druskin2021lippmann}.
The new system for $p$ is
\be\label{eq:oper1d}
\delta F= \int W(x) p(x) dx
\ee
where
$$\delta F=[(F_0-F)(\lambda_1),\ldots, (F_0-F)(\lambda_m), \frac{d}{d\lambda}(F_0-F)(\lambda_1),\ldots, \frac{d}{d\lambda}(F_0-F)(\lambda_m)],$$
and
$$W=[\uu {u}^0(\lambda_1),\ldots,
\uu {u}^0(\lambda_m),\frac{d}{d\lambda}(\uu{u}^0)|_{\lambda=\lambda_1},\ldots,
\frac{d}{d\lambda}(\uu {u}^0)|_{\lambda=\lambda_m}]$$ 
are $2m$-dimensional vectors of functions on $(0,L)$. Recall that $\uu$  is computed by (\ref{internal1d}) directly from the data without knowing $p$; thus, the nonlinear system (\ref{eq:int1d}) becomes linear. We refer to (\ref{eq:oper1d}) as a Lippmann-Schwinger-Lanczos system. 

\end{enumerate}

\section{Connection between frequency domain LSL and time snapshots}\label{ssec:time}
Consider the time domain wave problem corresponding to ${\cal L}$ from \eqref{eq:1D}
\[{\cal L} w(x,t)+w_{tt}(x,t)=0
\]
with the same boundary condition as \eqref{eq:1D} and initial 
conditions 
\[w|_{t=0}=g, \qquad w_t|_{t=0}=0.\]
Given regular enough $p$ and a localized initial pulse, $g(x)\approx \delta(x)$, we have that  
\[w\approx\delta(x-t)+\eth(x,t),\]
where  $\eth(x,t)$ is a smooth function and satisfies the causality principle, that is,  $\eth=0$ for $x\ge t$,  and $\eth\equiv 0$ if and only if $p\equiv 0$. 
Consider time snapshots $w(x,\tau i)$ for $i=1,\ldots,m$, for $\tau$ consistent with the effective width  of $g(x)$. Let us perform Gram-Schmidt orthogonalization on these snapshots in sequential order. This is the foundation of all data-driven inversion algorithms in the time domain: due to causality, the orthogonalized time snapshots are approximately $\delta(x-\tau i)$ \cite{druskin2016direct}, the same as those from the background medium. Therefore, the orthogonalized time snapshots depend only weakly on $\eth$, and likewise on $p$. Although the snapshots are unknown, the transformation that orthogonalizes them is obtained from the data generated Gramian.  

 In the frequency domain, the snapshots are not causal, so sequential Gram-Schmidt orthogonalization will not lead to weak dependence on $p$. In \cite{BoDrMaMoZa}, weak dependence on $p$ was extended to the frequency domain by instead using Lanczos orthogonalization. To see why this works, consider the approximation of the time snapshots $w$ in the Galerkin framework as
\[w(x,t)\approx V(x) d(t),\]
where the time dependent coefficients $d\in \RR^m$ satisfy
\begin{equation} \label{timeROM} S d(t)+ M d(t)_{tt}=0,
\qquad w|_{t=0}=b, \qquad w_t|_{t=0}=0,\end{equation}
which is a time-domain (the wave) variant of \eqref{eq:ROM}.
Thus, snapshots of its solution, \[ V d( \tau i )\approx\delta(x-\tau i)+\eth(x,\tau i),\]
similar to the full time domain discussed earlier. Their orthogonalization creates a basis of approximate delta functions, however, their peaks, $t_i$, are not equidistantly spaced. This is because the approximation properties of the basis $V$ deteriorate away from the origin. It was observed in \cite{BoDrMaMoZa} that each peak $t_i$ is located within $i$-th dual volume of the corresponding optimal finite-difference grid.  We note that it was also observed and proven for a model case in \cite{borcea2005continuum} that the corresponding optimal grids are only weakly dependent on the medium.

Since the solution to (\ref{timeROM}) is a cosine, $d(\tau i)$ satisfies {\it exactly} the second order finite-difference scheme
\be\label{eq:fd} d[\tau(i+1)]=(2I-\tau A)d[\tau i]-d[\tau(i-1)], i=i,\ldots, m-1,\ee 
\[d(0)= M^{-1} b, \quad d(\tau)=d(-\tau), \]
where $A= M^{-1}S$.
Thus, by construction, each $d(\tau i)$ is from the Krylov subspace and for each $i$,
 $$\hbox{span} \{ d(0), d(\tau),\ldots, d(\tau i)\} =\hbox{span} \{M^{-1}b, AM^{-1}b, \ldots,  A^iM^{-1}b\}. $$
 Hence, when we do the Lanczos algorithm as described in Step 2, we are performing exactly the sequential Gram-Schmidt orthogonalization of the ROM projected time snapshots.
 


\section{The regularized Lippmann-Schwinger-Lanczos method}\label{reg-LSL}
Data-driven formulas for mass and stiffness matrices (\ref{eq:massmtr},\ref{eq:stifmtr})  are only valid up to the error in the data. For ill-conditioned snapshot spaces, the small singular values of $V$ and ${\sqrt{L}}V$ may not be estimated well by the eigenvalues of $M$ and $S$. Furthermore, $M$ and $S$ may even become indefinite, which would lead to losing the Hermitian property of the matrix pencil $(S,M)$, fundamental to ROM imaging methods. To avoid the resulting instability, we regularize the problem by truncating the spectral decomposition of the Gramian matrix $M$ and projecting the problem into the subspace spanned by the dominant eigenvectors. Truncated SVD is an efficient regularization method that allows constructing a reduced order model using only eigenvectors corresponding to the dominant eigenvalues as the basis of the projection subspace.   
The Gramian truncation was used for inverse scattering in the time domain \cite{borcea2019robust} and for balanced truncation of data-driven  ROMs in the frequency domain \cite{gosea2022data}. The regularized algorithm is presented in section  \ref{Meigenvectors}. 
\subsection{ROM regularization via Gramian truncation}
 Let $Z$ be a matrix of the dominant eigenvectors of the mass matrix $M$:
\begin{equation} \label{Meigenvectors} 
Z= [z_1,\ldots, z_l]\in \RR^{m\times l},
\end{equation} 
with $z_k, k = 1, \ldots, l$,  for $ l \le m$, being the eigenvectors of $M$ corresponding to its positive eigenvalues 
$\sigma_k, \, 1\le k\le l$, that are greater than a threshold $\alpha$,  
  $\sigma_k \ge\alpha$. The parameter $\alpha$ 
is related to the level of noise in the data and determines the size of the ROM by the number $l$ of the basis vectors.

We will use the subspace spanned by these eigenvectors for constructing the reduced order model by projecting the problem into this subspace. Let 
\begin{equation}\label{eq:MSbt}
\tM=Z^*MZ= \tilV^*\tilV, \quad  \tS=Z^*SZ=\tilV^*L\tilV, \end{equation} \begin{equation}  \tb=Z^*b, \quad \mbox{and} \quad \tilV=V Z\in \RR^{\infty\times l}.
\end{equation}
Then, the Galerkin system (\ref{eq:ROM}) determining the truncated ROM can be written as 
\be\label{eq:ROMt}
(\tS+\lambda \tM)\tc=\tb
\ee
The Galerkin solution for the truncated problem can be expressed as 
\be\label{eq:ROMsolfr} 
u(\lambda)\approx \tilde{ u}(\lambda)= \tilV  \tc(\lambda)=\tilV\left(\tS+\lambda \tM\right)^{-1}\tb.
\ee
Due to the truncation, $\tilde{u} $ does not exactly match the data. We observe computationally that this misfit does not exceed the order of the measurement error, however, a more detailed error analysis is needed.

Since $\tM$ is the diagonal matrix of positive eigenvalues of $M$, this truncation guarantees the Hermitian property of the matrix pencil $(\tS, \tM)$ and allows to use the $\tM$-Hermitian  Lanczos algorithm  with positive-definite $\tM$. This algorithm is isomorphic to the standard real symmetric Lanczos method with matrix $\tM^{-1/2}\tS\tM^{-1/2}$ and  results in  
 the computation of the internal solution without the Lanczos algorithm breaking down. 

\subsection{Regularization of the background model} 
In principle, there is less need to regularize the ROM for the background model since it  can be computed exactly from the simulated data, and all the eigenvalues of the corresponding mass matrix must be positive definite up to the precision of the computation. However, for the LSL approach, we need it to be consistent with the data-generated ROM. Hence we need to project it on the same $Z$ as above (\ref{Meigenvectors}) 
to obtain $\tM _0,\tS_0\in\RR^{l\times l}$ and $\tb_0\in \RR^{l}$ where 
\begin{equation}\label{eq:MSbt0}
\tM_0=Z^*M_0Z= \tilV_0^*\tilV_0, \quad  \tS_0=Z^*S_0Z=\tilV_0^*L_0\tilV_0, \end{equation} \begin{equation}  \tb_0=Z^*b_0, \quad \mbox{and} \quad \tilV_0=V_0 Z\in \RR^{\infty\times l}.
\end{equation}
Then, performing $m$ steps of the  $\tM_0$-Hermitian  Lanczos algorithm corresponding to  matrix $\tilde{A}_0=\tM _0^{-1}\tS_0$ and initial vector  $\tM _0^{-1}\tb _0$, we compute the orthogonalized snapshot basis for the background
$$\tilV_0 Q_0= V_0ZQ_0.$$ 

\subsection{Internal solutions for the truncated problem}
As we discussed in subsection \ref{ssec:time}, the orthogonalized snapshot columns are weakly dependent on $p$.   That is, 
\[VZQ \approx V_0ZQ_0, \]
from which we have that 
\[ VZZ^*\approx V_0ZQ_0Q^{-1}Z^*.
  \]
 Assuming that the projector $ZZ^*$ has good approximation properties on $V$, i.e,  $VZZ^*\approx V$, this implies that 
\be\label{eq:snaps}  V  \approx V_0ZQ_0Q^{-1}Z^*.
\ee
This approximation is obtained directly from the data and used to construct the internal solutions.

\subsection{Inversion using Lippmann-Schwinger}
 We consider the Lippmann - Schwinger formulation in order to solve the inverse problem by using the internal snapshots described above. 
Recall that 
	\be  \label{eq:int}
	F_0(\lambda_j)-F(\lambda_j)=\int {u}^0_j(x ) u_j(x )p(x)dx, \\
	\qquad j=1,\ldots, m \ee
where $u_j$ and $u^0_j$ are the perturbed (unknown) and background solutions, respectively.
With help of (\ref{eq:snaps})  we transform  \eqref{eq:int} into
\be\label{eq:LSL}  
F_0(\lambda_j)-F(\lambda_j)\approx \int {u}^0_j(x )\uu_j p(x)dx, \\
\qquad j=1,\ldots, l ,\ee
 where  \[\uu_j = e_j^*V_0ZQ_0Q^{-1}Z^*\]  can be directly computed from the data without knowing $p $, thus transforming the nonlinear system (\ref{eq:int}) into the linear system (\ref{eq:LSL}).  We can write the system \eqref{eq:LSL} in operator form   as
 \be\label{eq:oper}
 \delta\bF=\langle W, p\rangle 
 \ee
 where 
 $$
 \delta \bF=[F_0(\lambda_1)-F(\lambda_1),\ldots, F_0(\lambda_l)-F(\lambda_l)]^*\in\RR^{l},
 $$
   and     
 $$
 W=[\uu_1(x)u^0_1(x), \ldots,   \uu_l(x) u^0_l(x)]
 $$ 
 is an $l$-dimensional vector valued function on $(0,L)$.
 \begin{remark} \label{rem:mpoints}
  In the original LSL method, both  $F(\lambda_j)$ and their derivatives were used to construct the ROM and solve  \eqref{eq:oper}. In the regularized formulation,  the derivatives are used for the ROM construction; however,  their addition as data in \eqref{eq:oper} is not necessary as this does not improve the quality of the solution. 
 This simplifies the regularized LSL formulation, further reducing the size of the problem. Indeed, instead of a system of size $2 m$ that is necessary to solve using the data driven LSL method, the presented regularized approach only requires solving a linear system of size $l, l<m$.  
 \end{remark}

\section{Regularized LSL for multidimensional MIMO problems}\label{MIMO-LSL}
We consider a partial differential equation in the domain $\Omega\in\RR^d$ for a scalar function $u$
\be\label{d-Schrod}
({\cal L}+\lambda I )u^{(r)} =g^{(r)} \ \ \ \mbox{in} \ \ \Omega,  \qquad  \frac{\partial u^{(r)}}{\partial \nu}\large|_{\partial \Omega}=0, \ee for   $r=1,\ldots, K,$ 
where 
\[{\cal L}=-\Delta  +p I,\]
$\nu$ is the normal to the boundary of $\Omega$, 
and the sources $g^{(r)}$ are localized and supported near or at an accessible part of the boundary $\partial \Omega$.
As in the SISO case, for simplicity, we assume that  $\lambda$ is real, extension to the complex case is straightforward \cite{druskin2021lippmann}.

Define the vector of source functions $$G=[g^{(1)}, g^{(2)}, \ldots, g^{(K)}]$$  and corresponding solutions 
$$U=[U^1, U^2, \ldots, U^K],$$ 
which are  understood as semi-infinite matrices or vectors of continuous source and solution functions, respectively, 
$U^r \in \mathbb{R}^{\infty\times K}$, 
 corresponding to a particular spectral value $\lambda_j,  j=1, \cdots, m$, and $U^r = U^r(\lambda_j)$. 
Then the multiple-input multiple-output (MIMO) transfer function is a $K\times K$ matrix-valued function of $\lambda$
	\be
	F(\lambda)= \langle G,U\rangle =\langle G,({\cal L}+\lambda I )^{-1} G\rangle 
	\in\RR^{K\times K},
	\ee
where  
$\langle, \rangle$  represents the continuous  $L^2(\Omega)$  matrix-valued inner product, 
$$
\langle G,U\rangle = \int_\Omega G^* U dx. 
$$  
We consider the inverse problem with data given by $2m$ real symmetric $K\times K$ matrices, i.e., for real $\lambda$ the data are 
$$
F(\lambda)|_{\lambda=\lambda_j}\in\RR^{K\times K}, \quad \mbox{ and} \qquad  
\frac{F(\lambda)}{d\lambda}|_{\lambda=\lambda_j} \in \RR^{K\times K}, \qquad  j=1, \cdots, m.  
$$

\subsection{Construction of the MIMO data-driven ROM}
 The Galerkin projected system for the MIMO data-driven Reduced Order Model is given by the block analog of (\ref{eq:ROM})
 \be\label{eq:ROMMIMO}
(S+\lambda M)C = B.
\ee
Here $S,M\in\RR^{mK\times mK}$ are Hermitian positive definite matrices,  $B\in\RR^{mK \times K}$,  and $C=C(\lambda)\in\RR^{mK\times K} $ is a matrix-valued function of $\lambda$, corresponding to coefficients of the solutions with respect to the basis $V$ of exact solutions. Stiffness and mass matrices 
consist of $m \times m$ blocks $S_{ij}$ and  $M_{ij}$ of size $K\times K$  given by 
\be  \label{blockSM}
 S=(S_{ij}=\langle \nabla {U}_i,\nabla U_j\rangle +\langle p{U}_i,U_j\rangle),  
\qquad
  M=(M_{ij}=\langle {U}_i,U_j\rangle).
\ee  
They can be computed directly from the data by the block versions of (\ref{eq:massmtr}), (\ref{eq:stifmtr}). 
 Similar to the SISO case, the Galerkin projected system (\ref{eq:ROMMIMO}) results in the reduced order model and the ROM transfer function  \[ \tF(\lambda)=B^*C(\lambda) \] which 
matches the data exactly. 
The truncation leads to truncated ROM 
\be\label{eq:ROMtMIMO}
(\tS+\lambda \tM)\tC=\tB
\ee
where 
\be\label{eq:MSbtMIMO}
\tM=Z^*MZ= \tilV^*\tilV, \ \tS=Z^*SZ=\tilV^*L\tilV, \ \tB=Z^*B, \quad \tilV= V Z\in \RR^{\infty\times l}
\ee
are the MIMO analogues of (\ref{eq:MSbt}) and (\ref{eq:ROMt}), and $l$ is the dimension of the reduced problem,  $l \le mK$.

\subsection{Data-generated internal solutions for the MIMO case}
The data-generated solutions for the MIMO case are constructed following similar steps as in the SISO case but performed block-wise.  The localization properties for this case were investigated in \cite{druskin2018nonlinear}. 
It can be performed via $m$ steps of the {$M$-Hermitian}  block-Lanczos algorithm with matrix $$A=\tM^{-1}\tS$$ and initial block  vector $\tM^{-1}\tB$.

\begin{remark} The orthogonalization in the block-Lanczos algorithm is not unique; the blocks can be constructed in different ways. Here we use the version with a polar decomposition within the blocks  \cite{druskin2018nonlinear}.  
\end{remark}  

From the orthogonalization we obtain the  $M$-orthonormal Lanczos block-vectors $q_i\in \RR^{mK\times K}$ which are the  block counterparts of (\ref{eq:lancz}) 
	$$ Q=[q_1, q_2,\ldots, q_m]\in{\RR^{mK\times mK}},$$ where $$ q_1=\tM^{-1}\tB(\tB^*\tM^{-1}\tB)^{-1/2}.$$
From this, we obtain the block-equivalent of \eqref{eq:snaps}
\be\label{eq:snapsB}  V  \approx V_0ZQ_0Q^{-1}Z^*,
\ee
which gives us the data-generated matrix of internal snapshots: 
\begin{equation}\label{eq:mimodatainternal} \mathbf{U} = V_0ZQ_0Q^{-1}Z^*.\end{equation}
	
\subsection{ MIMO inverse problem}
From the MIMO variant of (\ref{eq:transfer}), we obtain the block Lippmann-Schwinger equation \be\label{eq:intMIMO}F_0(\lambda_j)-F(\lambda_j)=\int U^*_0(x,\lambda_j ) p(x)U(x,\lambda_j )dx, \qquad j=1,\ldots, m,\ee
where the subscript $0$ corresponds to background solutions with $p=0$. Similar to the SISO case, precomputing $\mathbf{U}$ via the data-driven \eqref{eq:mimodatainternal}  results in the  linear system for $p$:
\begin{eqnarray} \label{eq:intdMIMO2}(F_0-F)|_{\lambda=\lambda_j} &=& \int U^*_0(x,\lambda_j ) p(x) \mathbf{U}(x,\lambda_j )dx. \end{eqnarray}
We note that the regularization allows us to use more spectral points and considerably reduce the size of the linear system in (\ref{eq:intMIMO}) at this step of LSL by not using the derivative data.  Although the derivative data was necessary to construct the ROM, adding it at this step does not improve the quality of reconstructions. 

\section{Extension to Helmholtz problem}\label{H-LSL} 


Helmholtz equation governs steady-state oscillations of the fields in homogeneous or non-homogeneous media; these oscillations could be electromagnetic, mechanical, or acoustic and arise in ultrasound, acoustics, electromagnetics, seismology, and other applications. 
As a model, consider the 2D Helmholtz problem 
\begin{equation} \label{helmholtz} -\Delta u^{(r)}(x) + \lambda n(x) u^{(r)}(x)  = g^{(r)} \ \ \mbox{in} \ \ \Omega\subset\RR^2, \ \ \ \ \ \frac{\partial u^{(r)}}{\partial \nu}|_{\partial \Omega =0 } \end{equation}
where the sources $g^{(r)}$ for $r=1,\ldots K$ are localized distributions supported near or at an accessible part of $\partial \Omega$. 

 Equation (\ref{helmholtz}) can be viewed as a quasistatic diffusive Maxwell equation or Laplace-transformed diffusion equation for the low-frequency electric field in a conducting medium; the data may correspond to $H_2$ optimal or another optimal approximation of the time-domain problem, e.g.,  \cite{druskin2013solution,druskin2009solution}. In this case, (\ref{helmholtz}) is sometimes called diffusion approximation,  
$u$ represents the electric field polarized normally to the plane $\RR^2$, and $n$ is the product of constant magnetic permeability and spatially varying conductivity. 
 Equation (\ref{helmholtz}) can also be viewed as Laplace transformed wave equation (with $\lambda$ being the square of the Laplace frequency), then $n(x)$ is a variable index of refraction or the inverse square of the variable wave speed. The latter analogy will be useful in the explanation of the numerical results.

The MIMO regularized LSL algorithm for this problem follows similar steps as for the Schrodinger LSL problem (\ref{d-Schrod}), with two important distinctions: 
\begin{enumerate}
    \item  The mass and stiffness matrices are computed by the same formulas (\ref{eq:massmtr}) and (\ref{eq:stifmtr}) directly from the data; the algorithm does not change here.  However, constructing variational Galerkin approximation, we see that the unknown coefficient $n(x)$  appears in $M$ instead of $S$ (compare with (\ref{blockSM})):
\be \label{HblockSM}
 S=(S_{ij}=\langle \nabla {U}_i,\nabla U_j\rangle+\langle {U}_i,U_j\rangle)
\qquad 
M=(M_{ij}=\langle n\, {U}_i,U_j\rangle).
\ee 
At the Lanczos orthogonalization step, the M-symmetric Lanczos algorithm is applied to matrix $A = M^{-1}S$ and initial vector $M^{-1}b$ to produce tridiagonal matrix $T$ and $M$-orthonormal Lanczos vectors $q_i$. 
However, in the Helmholtz case, the basis functions need to be orthogonalized with respect to the inner product with unknown weight $n$, see (\ref{HblockSM}), different from the background inner product, where we assume $n_0 =1$. Hence, we expect to see some errors that result in larger errors in the data generated internal solutions.

    \item In the case of the Helmholtz problem, the Lippmann-Schwinger equation is slightly different from the LS equation for the Schrödinger equation (\ref{eq:intdMIMO2}); to be consistent with the Helmholtz model, the data-driven Lippmann - Schwinger-Lanczos equation should be changed to: 
    \begin{eqnarray} \label{eq:intdMIMO2helmholtz}(F_0-F)|_{\lambda_j} &=& \lambda_j\int U^*_0(x,\lambda_j ) (n(x)-1)\mathbf{U}(x,\lambda_j )dx, \end{eqnarray}
with $\lambda_j$ present in front of the integral. 
 In the next section, we demonstrate reconstruction examples for the Helmholtz equation using the regularized LSL algorithm.    
\end{enumerate}

\section{Numerical Results}\label{Num}

In this section, we demonstrate the use of the Regularized Lippmann-Schwinger-Lanczos method (denoted Reg-LSL) for coefficient reconstruction.
We show that the improved stability of the Regularized Lippmann-Schwinger-Lanczos method allows us to use larger data sets and 
reconstruct the unknown potential in the presence of noise. We also compare its performance for both the Schrodinger and Helmholtz problems.

\subsection{Reconstructions in one dimension}

We first apply Reg-LSL to the SISO problem for both the Schrodinger and Helmholtz problems (the Helmholtz equation with positive $\lambda$ modeling diffusive electromagnetic problem). We compare the bases after Lanczos orthogonalization to those from the background, the corresponding internal solutions,  and the reconstructions of the perturbation.
Both problems are discretized on the interval $x\in[0,1]$ with a background perturbation of $p_0(x)=0$ for the Schrodinger problem and $n_0(x)=1$ in the Helmholtz problem.
The internal perturbation follows a Gaussian distribution scaled by a constant factor $\gamma=.125$, centered at $\mu=0.2$ and with standard deviation $\sigma=0.05$.
The Dirichlet data is collected at a single point source located at the origin for $\lambda = \{2,4,8,16,32,48\}$ in the Schrodinger case and $\lambda=\{2,4,8,16,32,48,64,96\}$ in the Helmholtz case. The synthetic data is computed by solving the forward problem using a centered finite difference scheme with step size $h=0.002$.
The regularization is enforced in both problems by truncating the singular values of the mass matrix for $\lambda_M\leq5e$-$12$. 
For the linear inverse problem in \eqref{eq:LSL}, we use the Moore-Penrose pseudo inverse with a truncation level of $6e$-$5$ for both the Schrödinger and Helmholtz problems. 

\begin{figure}
    \centering
    \begin{tabular}{cc}
    \includegraphics[width=.45\textwidth]{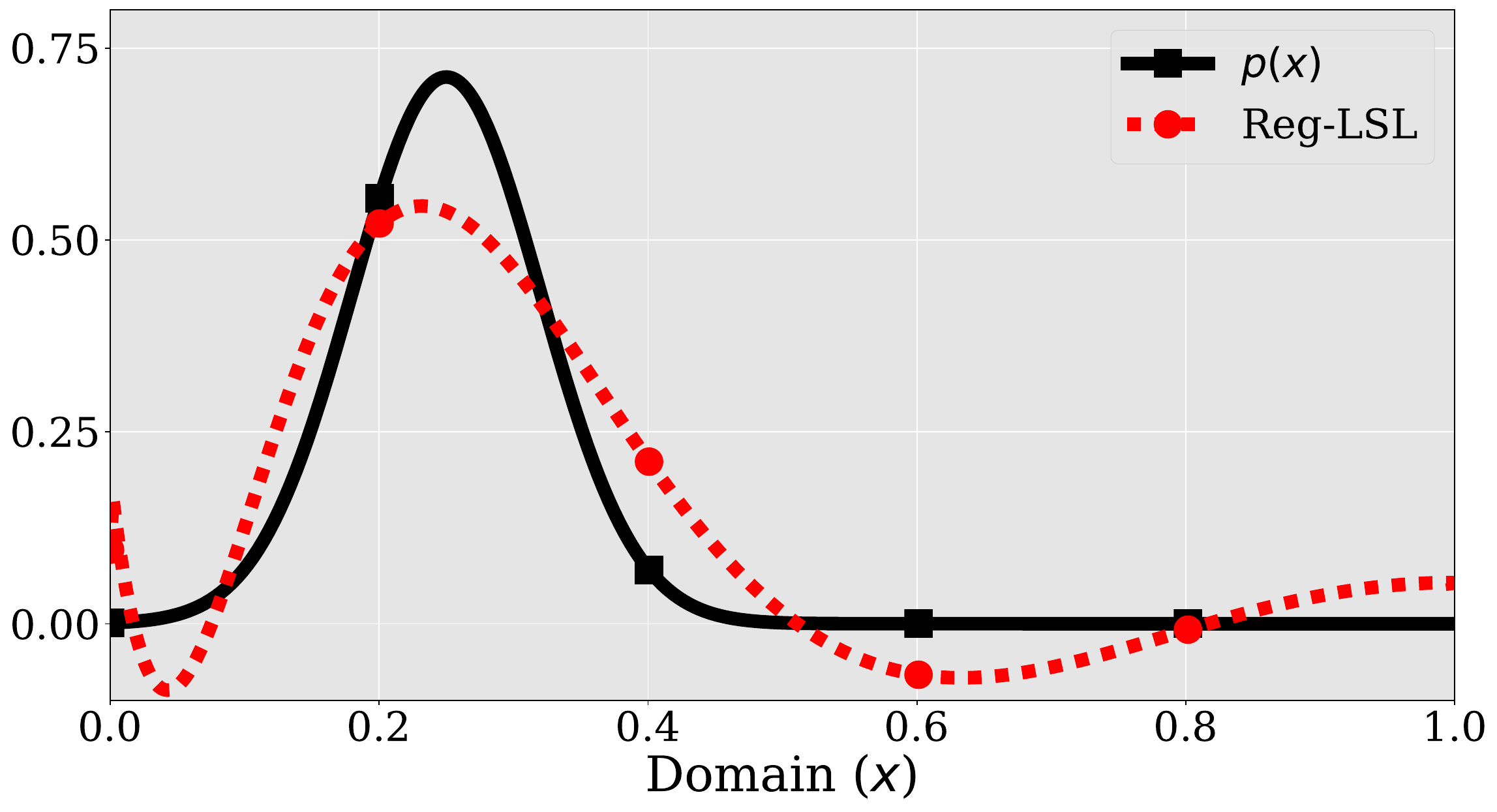} &
    \includegraphics[width=.45\textwidth]{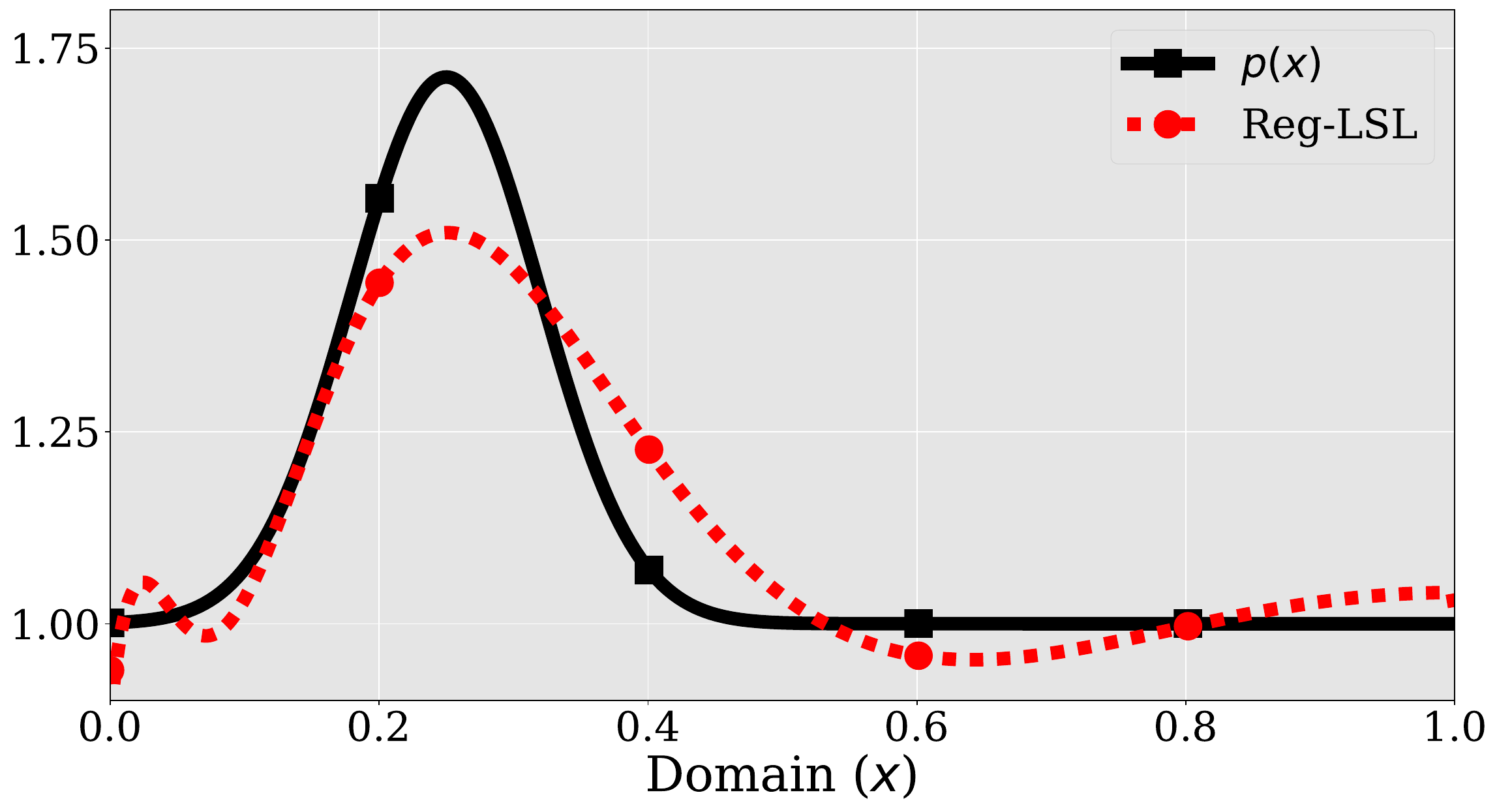} \\
    a) Schrödinger Problem & b) Helmholtz Problem \\
    \end{tabular}
    \caption{Perturbation (black-solid) and the RegLSL reconstructed perturbation (red-dashed) for a) the Schrödinger problem and b) the Helmholtz problem.
    }
    \label{fig:SISO:pert}
\end{figure}

Figure~\ref{fig:SISO:pert} illustrates the perturbation and reconstructed solutions for the a) Schrödinger problem and the b) Helmholtz problem. In both problems, we observe excellent reconstructions of the perturbation with few oscillations near the source. The Helmholtz reconstruction is more accurate outside of the support of the perturbation, at which point it begins to deviate due to the difference in conductivity, or equivalently, in the wave speed, affecting the slowness metric.

\begin{figure}
    \centering
    \begin{tabular}{cc}
    \includegraphics[width=.45\textwidth]{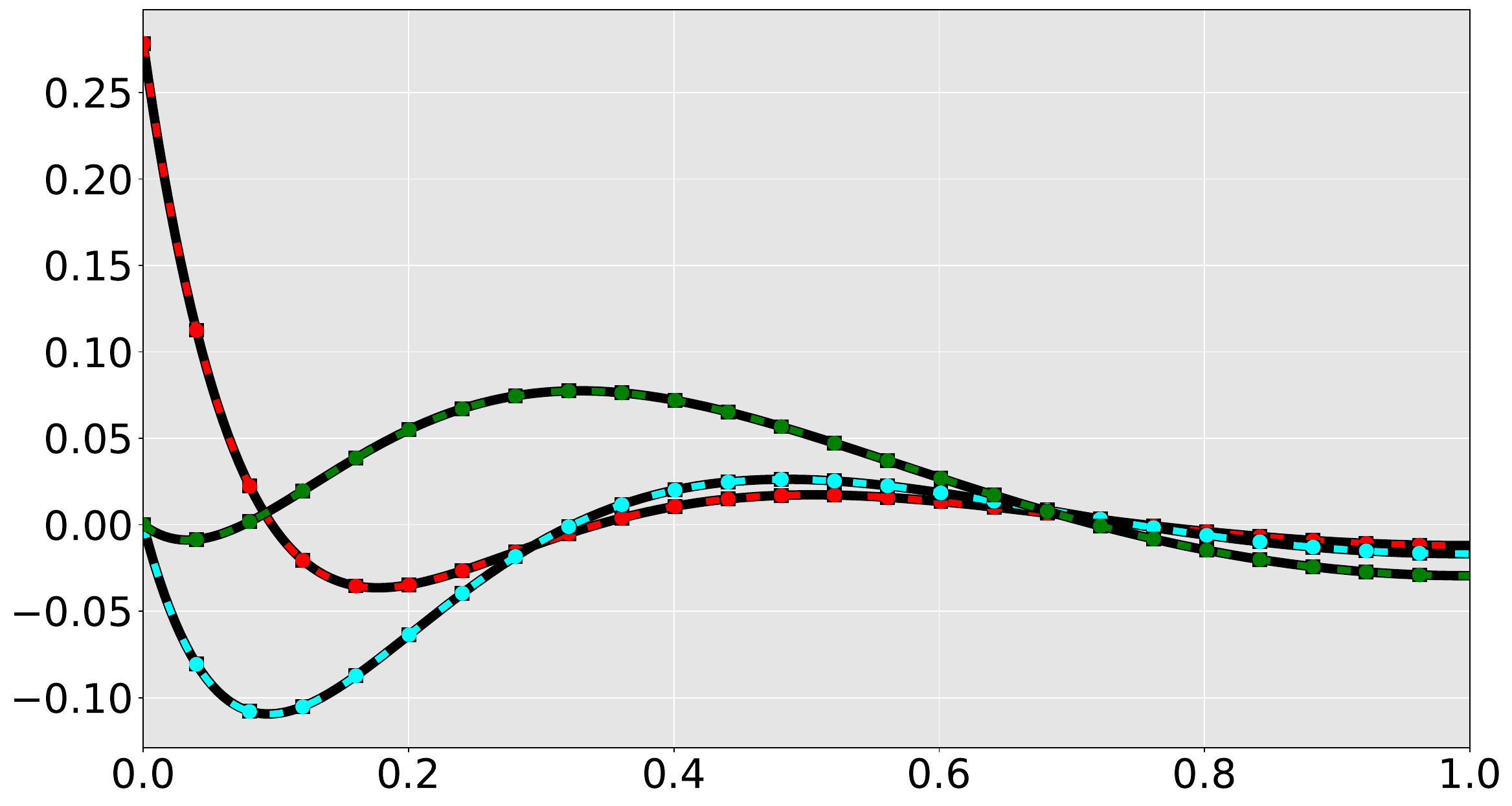} &
    \includegraphics[width=.45\textwidth]{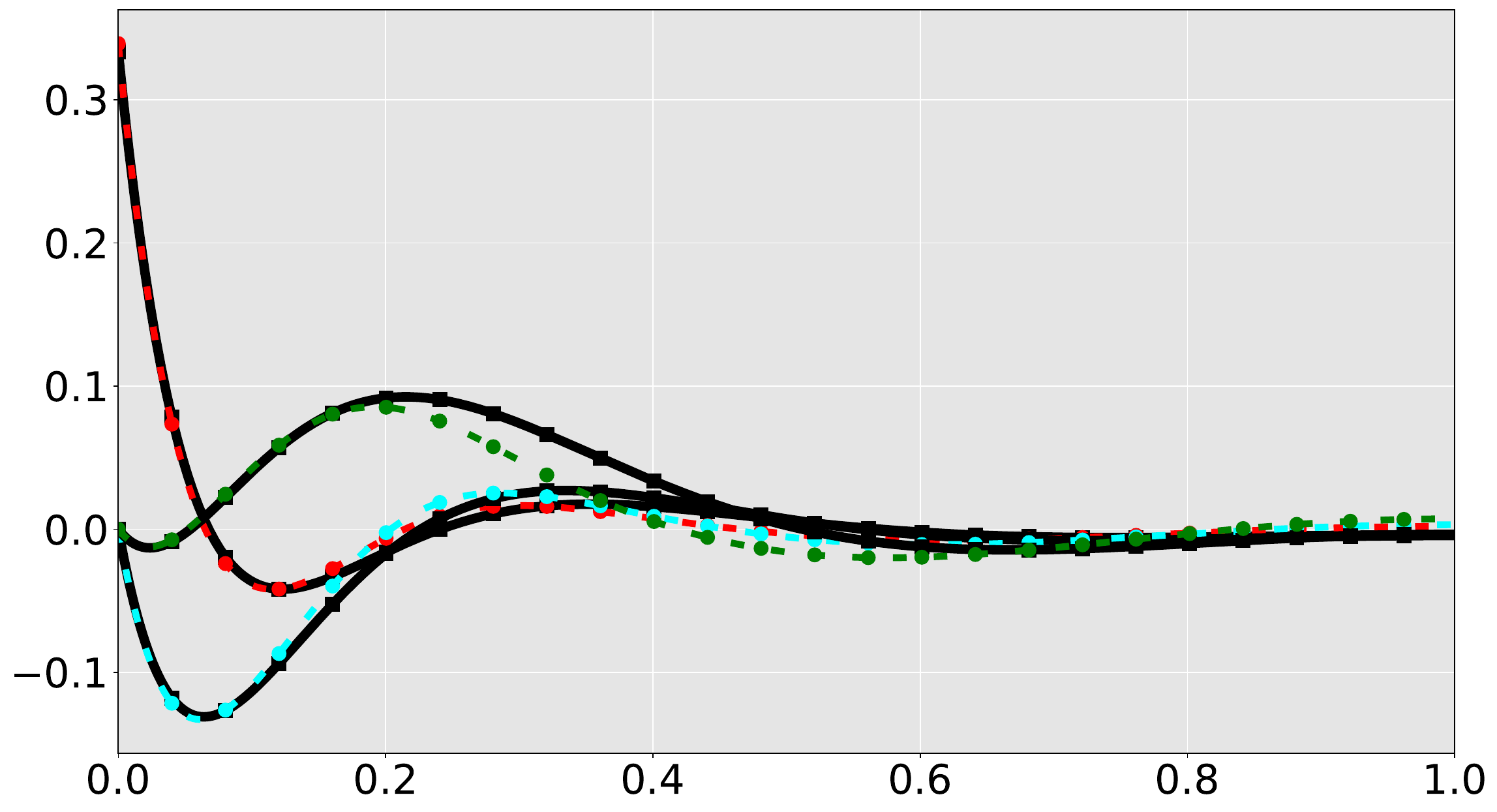} \\
    a) Schrödinger Problem & b) Helmholtz Problem \\
    \end{tabular}
    \caption{Comparison of three internal basis functions (black-solid) and their reconstructed counterparts (colored-dashed).}
    \label{fig:SISO:basis}
\end{figure}

The Lanczos orthogonalized basis functions for the SISO problems are shown in Figure~\ref{fig:SISO:basis}. We observe that the basis 
functions for the perturbed Schrodinger problem are almost exactly the same as those for the background problem; however, the basis functions for the Helmholtz equation deviate from the background basis in the area of the perturbation, again, due to the difference in the travel time metric.

\begin{figure}
    \centering
    \begin{tabular}{cc}
    \includegraphics[width=.45\textwidth]{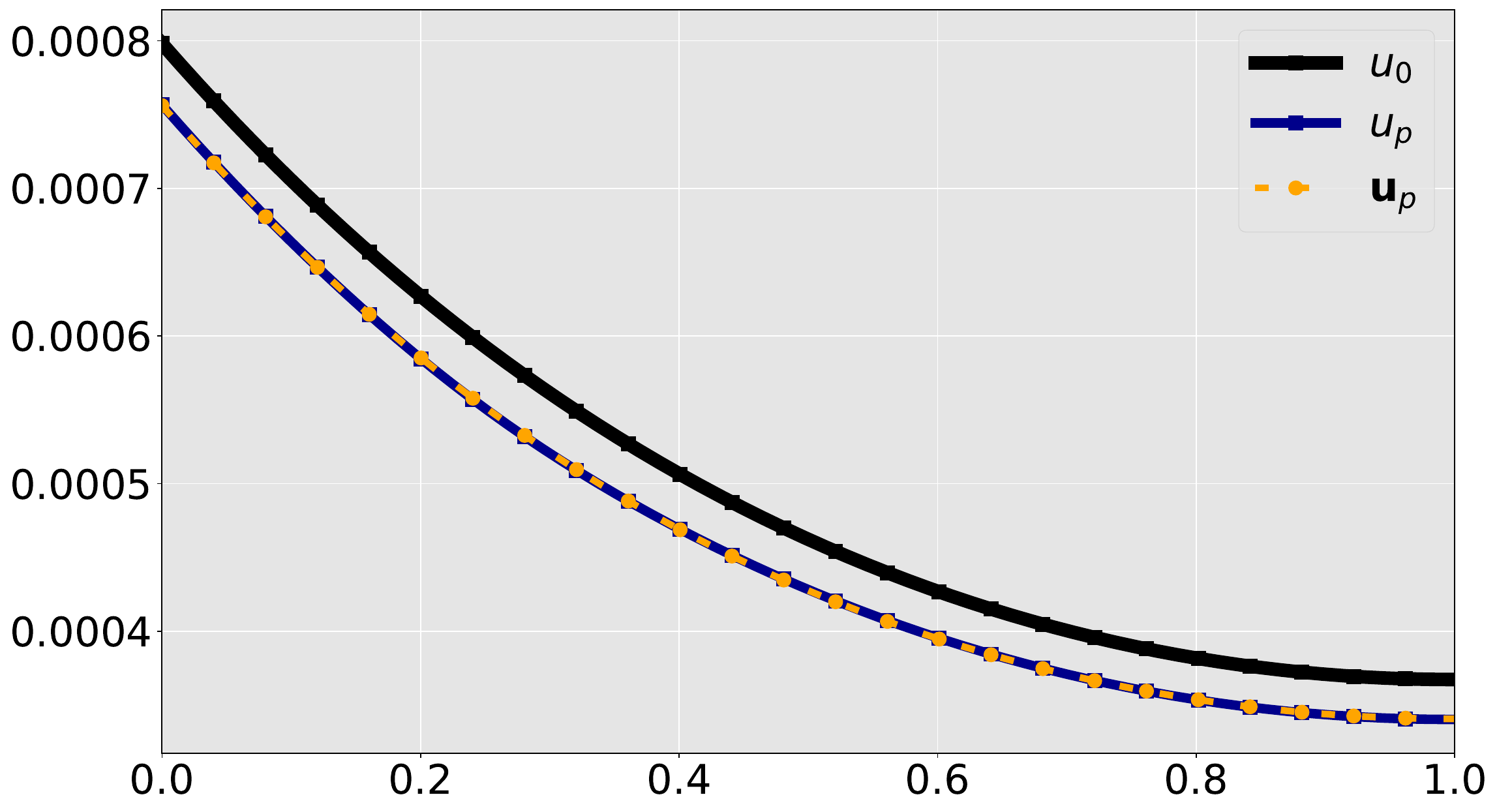} &
    \includegraphics[width=.45\textwidth]{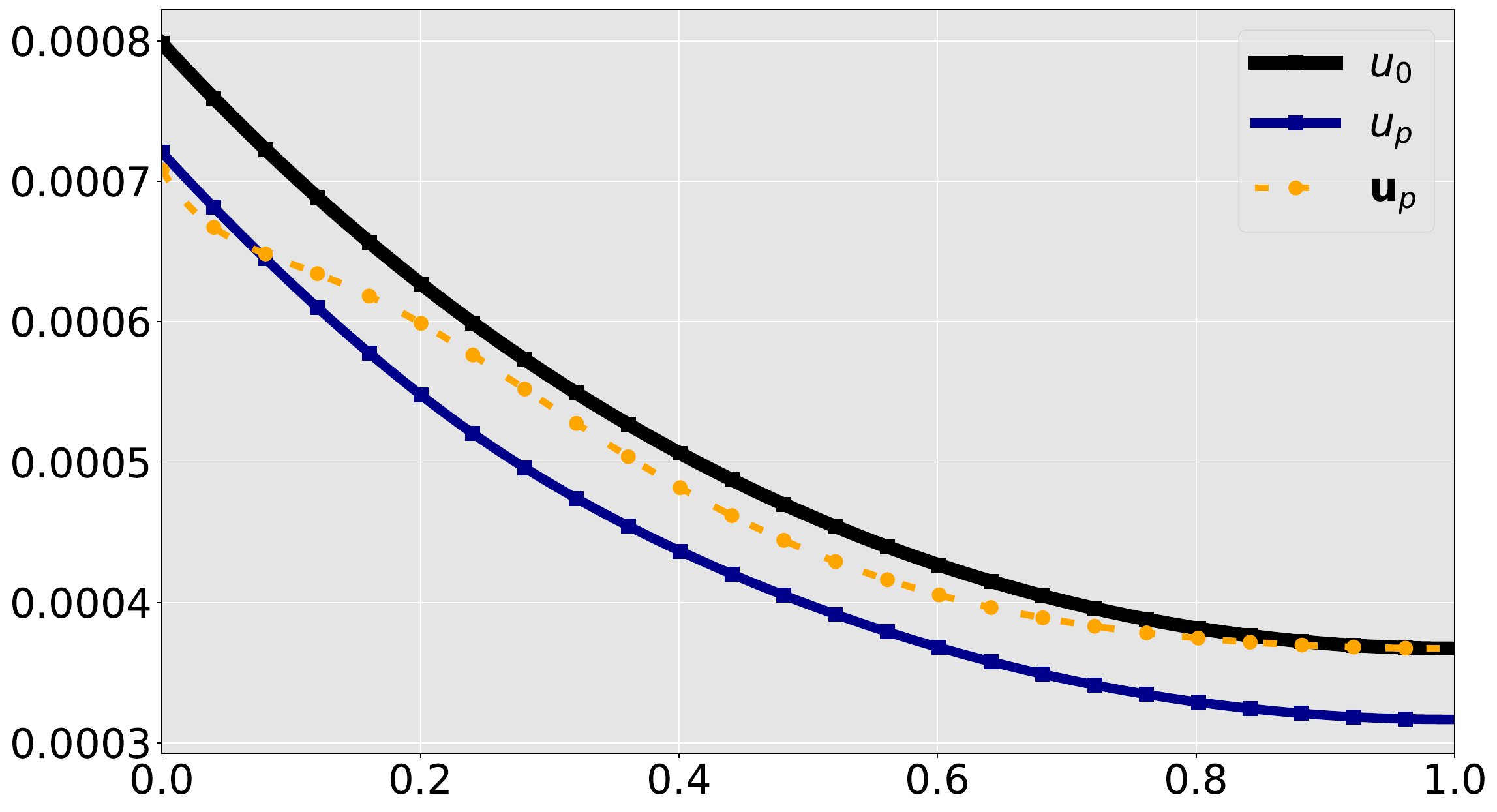} \\
    a) Schrödinger Problem & b) Helmholtz Problem \\
    \end{tabular}
    \caption{The background solution (black-solid) and the internal solution (blue-solid) with its reconstruction (yellow-dashed).}
    \label{fig:SISO:recon}
\end{figure}

Figure~\ref{fig:SISO:recon} shows the reconstruction of the internal solutions from data using the orthogonalized basis functions. 
The data generated solution for the Schrödinger problem shown in Figure~\ref{fig:SISO:recon}a, is nearly exact. In Figure~\ref{fig:SISO:recon}b, the internal solution deviates from the true one at the perturbation. Travel time is not affected in the interval between the source and perturbation, which is why we have an exact reconstruction of the internal solution there.  
 Similarly, a good accuracy of the reconstructions is achieved in the multidimensional setting (see below)  with localized conductive perturbations (corresponding to lower wave speed) since they do not affect travel time outside of the perturbation. Low conductivity (high wave-speed) perturbations, however, may work as waveguides, affecting the metrics globally. For this reason, direct LSL would not be applicable for such problems, for which one has to implement iterations \cite{doi:10.1137/22M1517342}

\subsection{Reconstructions for the MIMO problem}
In this section, we compare reconstructions for the MIMO Schrödinger and Helmholtz problems using Born, LSL, and Reg-LSL.
The MIMO problems presented in this work were solved using centered finite difference schemes with symmetrized Laplacian on the domain $(x,y)\in[-1,1]\times[-1,1]$ using a grid step $h_x=h_y=0.04$.
The Dirichlet data is collected at eight boundary source locations (two at each boundary) at evenly spaced intervals.

\subsubsection{Schrödinger Equation}

The Schrödinger problem has a background $p_0 = 0.0$ and a perturbation $p$ consisting of two Gaussian distributions.
The first perturbation is centered at $\mu_1=(0.2,0.5)$ with deviation $\sigma_1=(0.26,0.25)$, and the second is centered at $\mu_2=(-0.3,-0.5)$ with $\sigma_2=(0.2,0.18)$.
For Reg-LSL, the level of regularization is $5e$-$14$ for the Grammian truncation and $5e$-$4$ for the LS equation.
For Born, the regularization level is $7e$-$4$,  while the regularization level is $5e$-$3$ for LSL.  The data consists of the frequency set $\lambda=\{2,4,6,8,16,32,48\}$.

\begin{figure}[h!]
    \centering
    \begin{tabular}{cccc}
        \includegraphics[width=.21\textwidth]{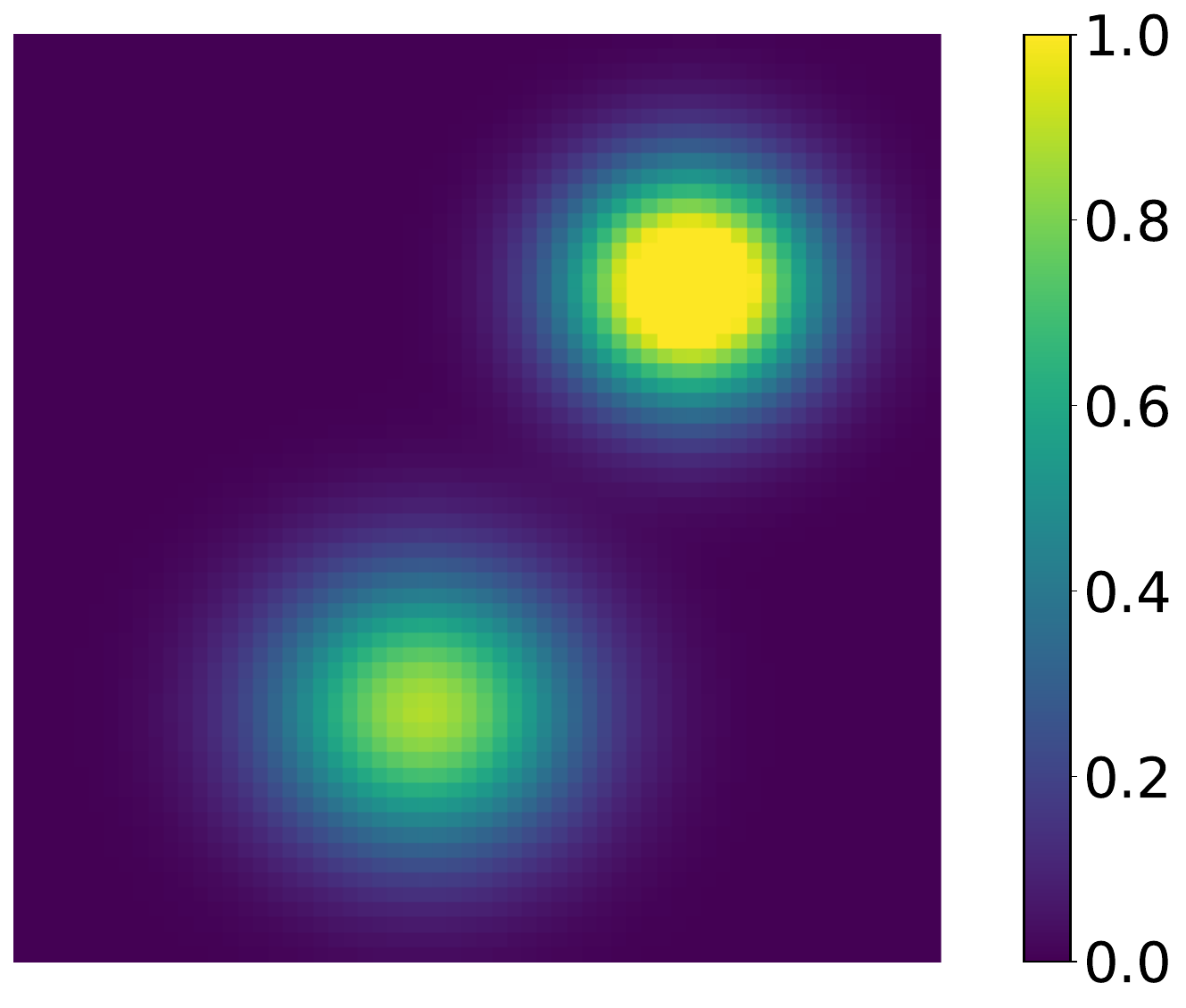} &
        \includegraphics[width=.21\textwidth]{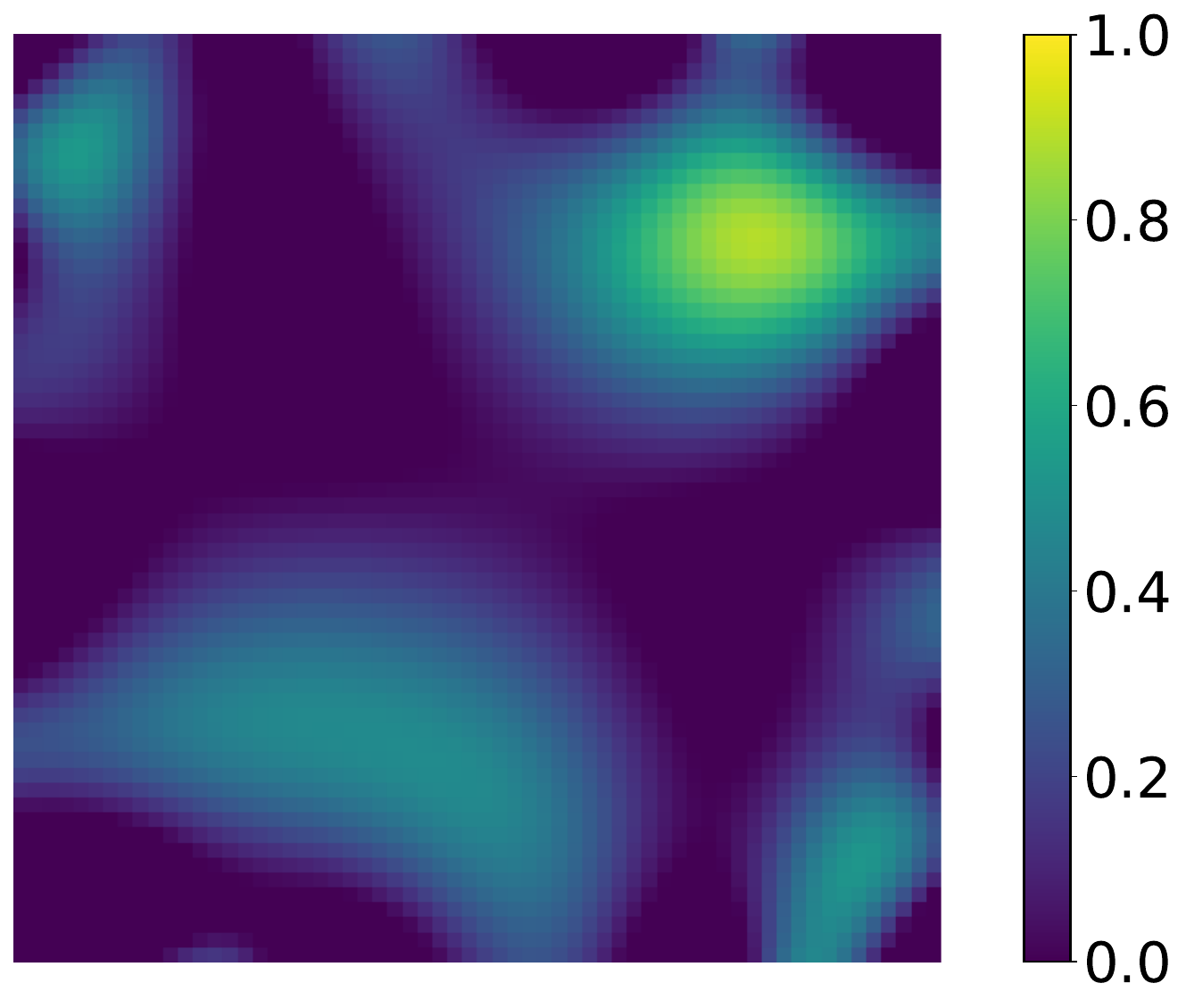} &
        \includegraphics[width=.21\textwidth]{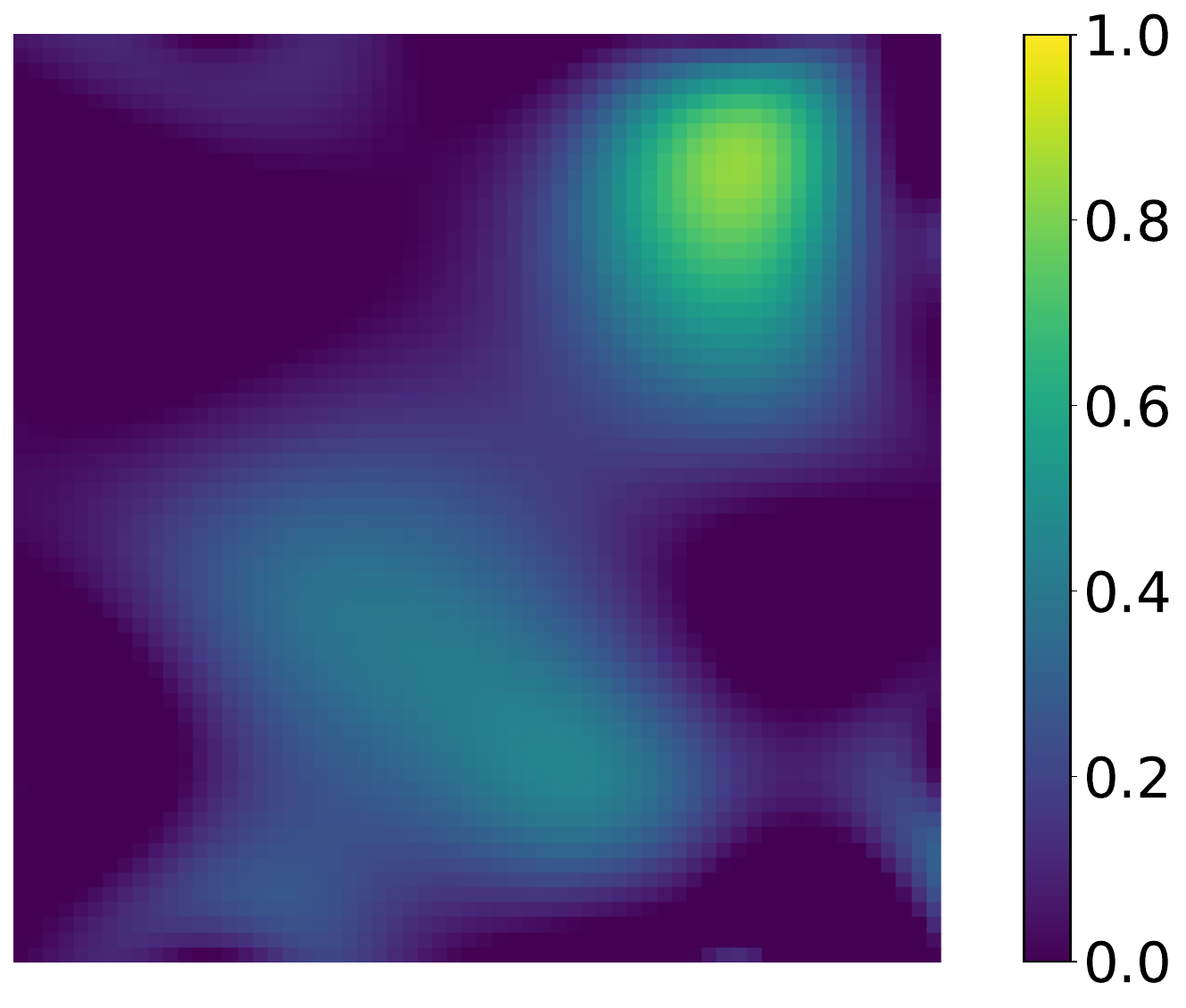} &
        \includegraphics[width=.21\textwidth]{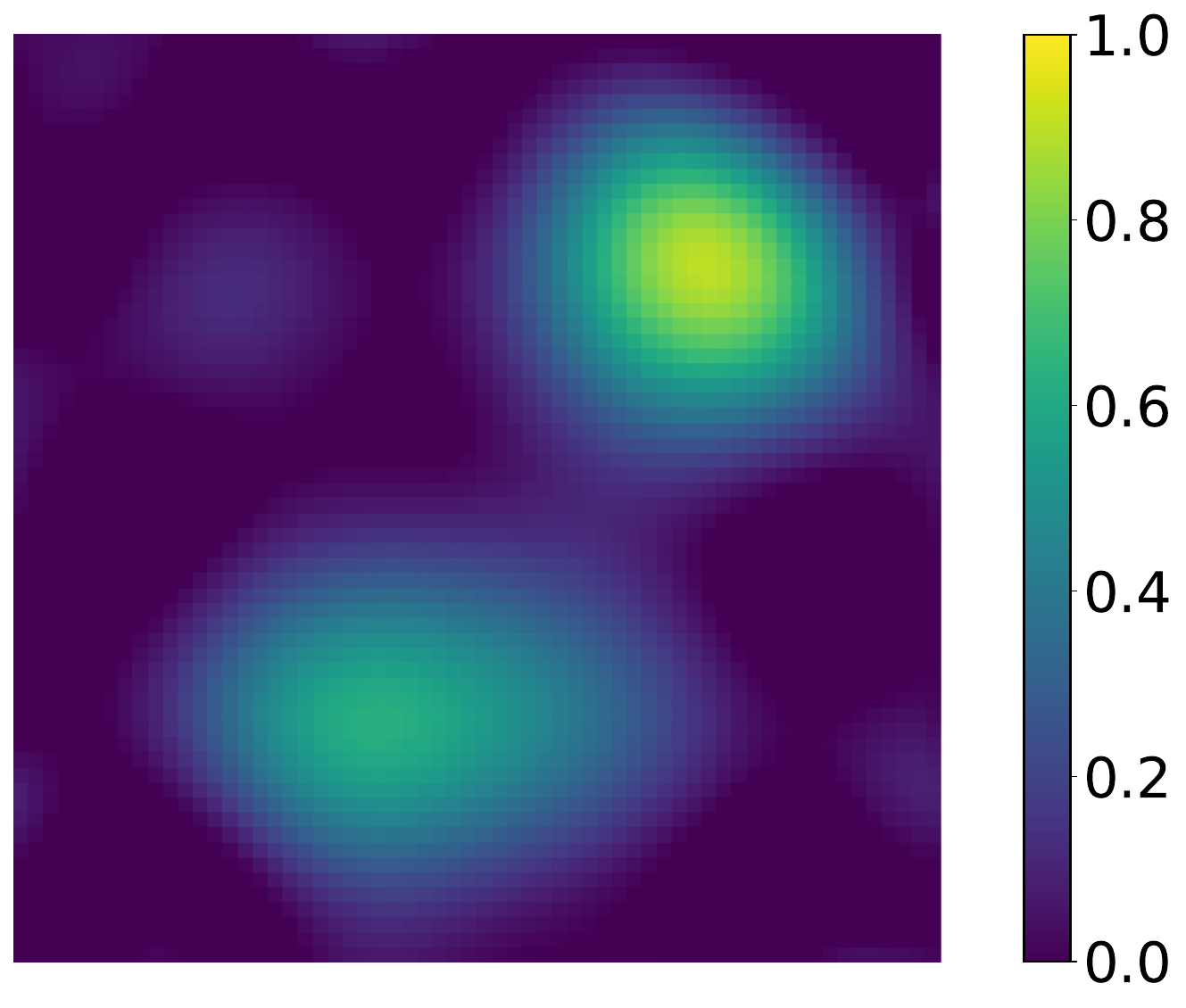} \\
        (a) True & (b) Born & (c) LSL & (d) Regularized LSL
    \end{tabular}
    \caption{Perturbation reconstructions for the MIMO Schrödinger problem.}
    \label{fig:mimo:schrodinger}
\end{figure}

Figure~\ref{fig:mimo:schrodinger} shows the reconstructions of the perturbations using each of the algorithms.
We observe that since this is a large number of frequencies for LSL, the standard LSL algorithm has reduced reconstruction quality.
Additionally, since this is a moderate contrast problem, the Born reconstruction has significant artifacts. 
By contrast, Reg-LSL has the lowest level of artifacts and has the highest fidelity to the true perturbation. 

\subsubsection{Helmholtz problem} 

For the quasi-stationary inverse conductivity problem modeled by the Helmholtz equation, we use the background $n_0=1.0$, with a perturbation consisting of three Gaussian distributions.
The first distribution is centered at $\mu_1=(-0.4,0.5)$ with standard deviation $\sigma_1=(0.16,0.15)$, the second is centered at $\mu_2=(-0.3,-0.4)$ with standard deviation $\sigma_2=(0.2,0.18)$, and the third is centered at $\mu_3=(0.4,0.2)$ with standard deviation $\sigma_3=(0.2,0.18)$.
For Reg-LSL, the Gramian truncation level is $1e$-$16$, and the level of regularization for its LS equation is $3e-4$.
The regularization levels for inverse Born and standard LSL are  $2e$-$3$ and $1e$-$3$, respectively. 

\begin{figure}[h!]
    \centering
    \begin{tabular}{cccc}
        \includegraphics[width=.21\textwidth]{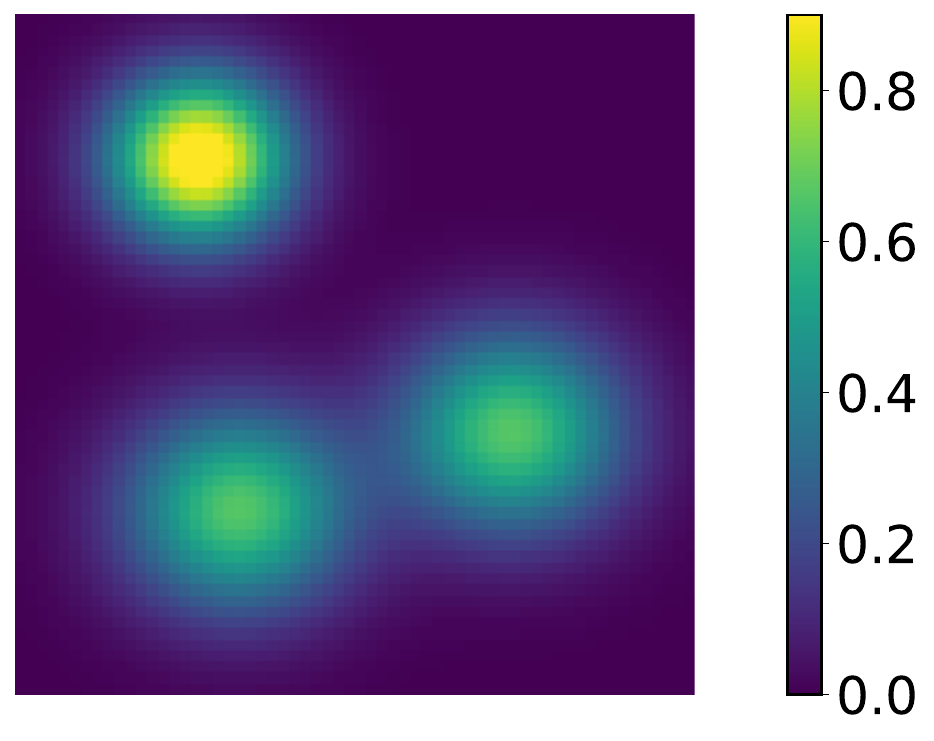} &
        \includegraphics[width=.21\textwidth]{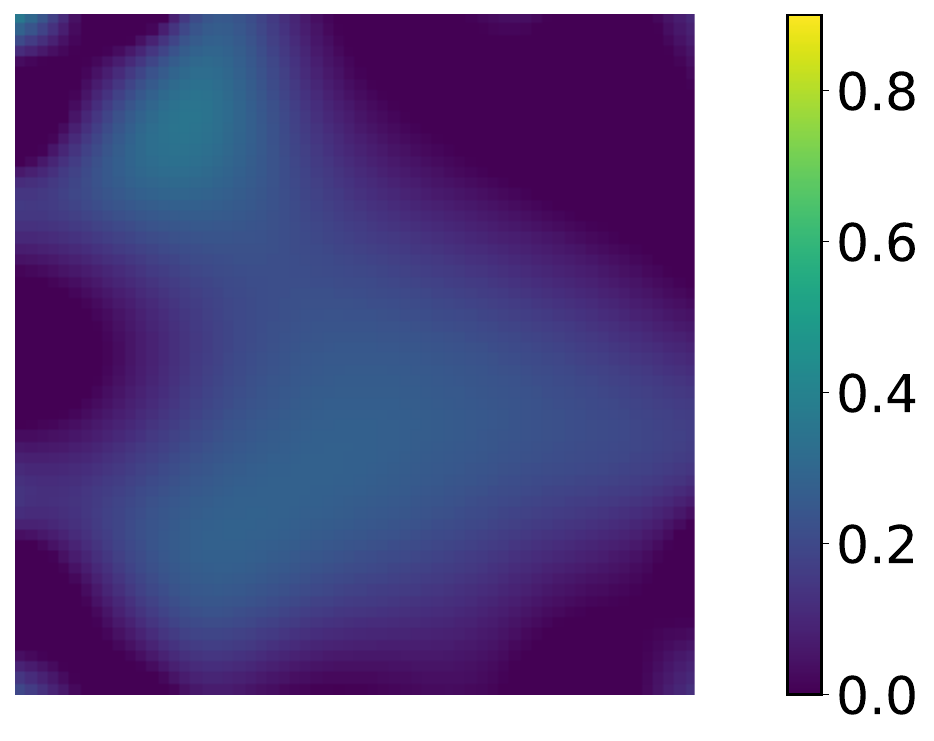} &
        \includegraphics[width=.21\textwidth]{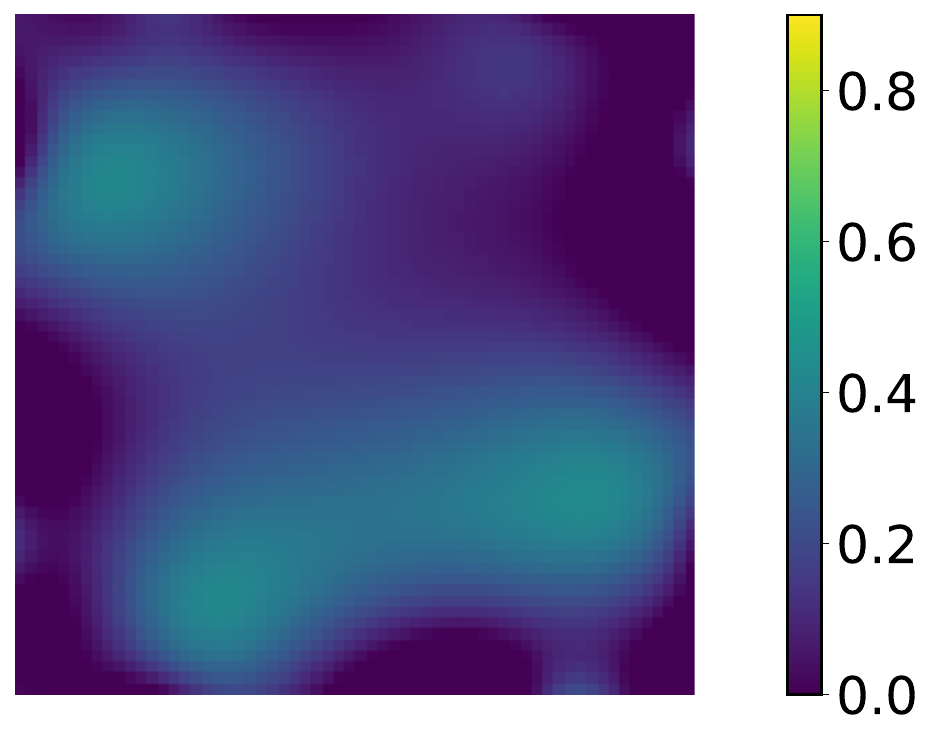} &
        \includegraphics[width=.21\textwidth]{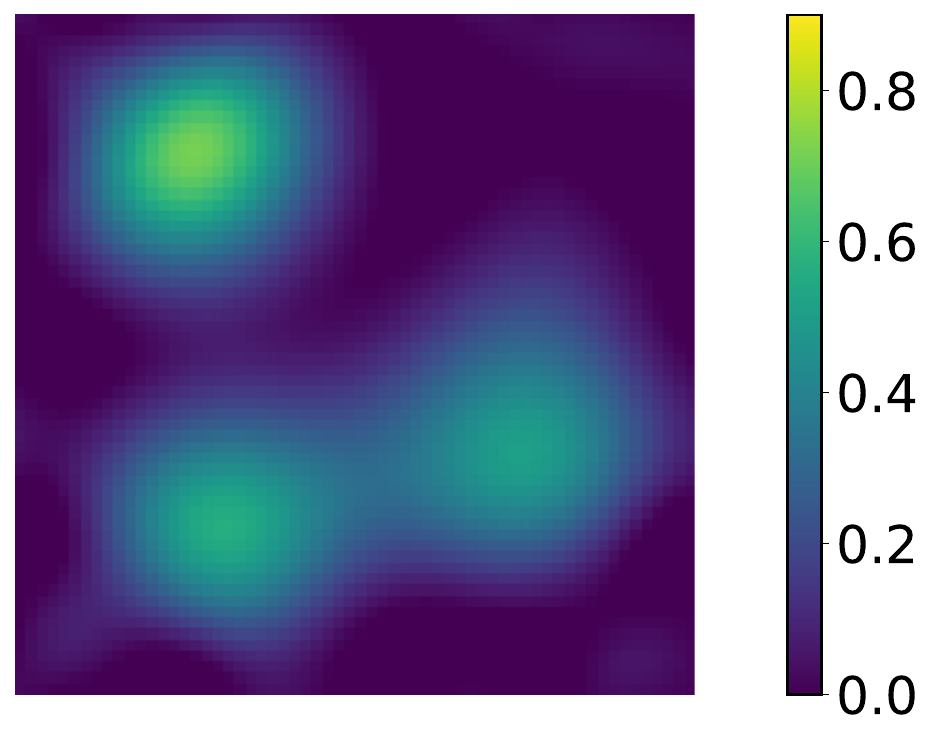} \\
        (a) True & (b) Born & (c) LSL & (d) Regularized LSL
    \end{tabular}
    \caption{Perturbation reconstructions for the MIMO Helmholtz problem.}
    \label{fig:mimo:helmholtz}
\end{figure}

Figure~\ref{fig:mimo:helmholtz} shows the reconstructions of the perturbation in $n$ using each of the algorithms.
The Born reconstruction has the lowest quality, followed by the standard LSL reconstruction. 
By contrast, Reg-LSL achieves the best reconstruction of the true perturbation. 

\subsection{Effects of regularization} In this section, we analyze the effects of the regularization level on the quality of the reconstructions.
Using the same MIMO problem configurations, we consider various thresholds for the regularization of the mass matrix.
For both the Schrödinger and Helmholtz problems, we used mass matrix regularization thresholds of $1e$-$4$,$1e$-$8$,$1e$-$12$, and $1e$-$16$.
For each of these thresholds, the levels of truncation for the linear part of the Schrödinger are $5e$-$4$, $5e$-$4$, $5e$-$3$, $2e$-$3$ respectively.
The levels of truncation for the linear inverse of the Helmholtz problem are $1e$-$3$,$1e$-$3$,$3e$-$4$, and $3e$-$4$, respectively.


\begin{figure}[h!]
    \centering
    \begin{tabular}{cccc}
        \includegraphics[width=.21\textwidth]{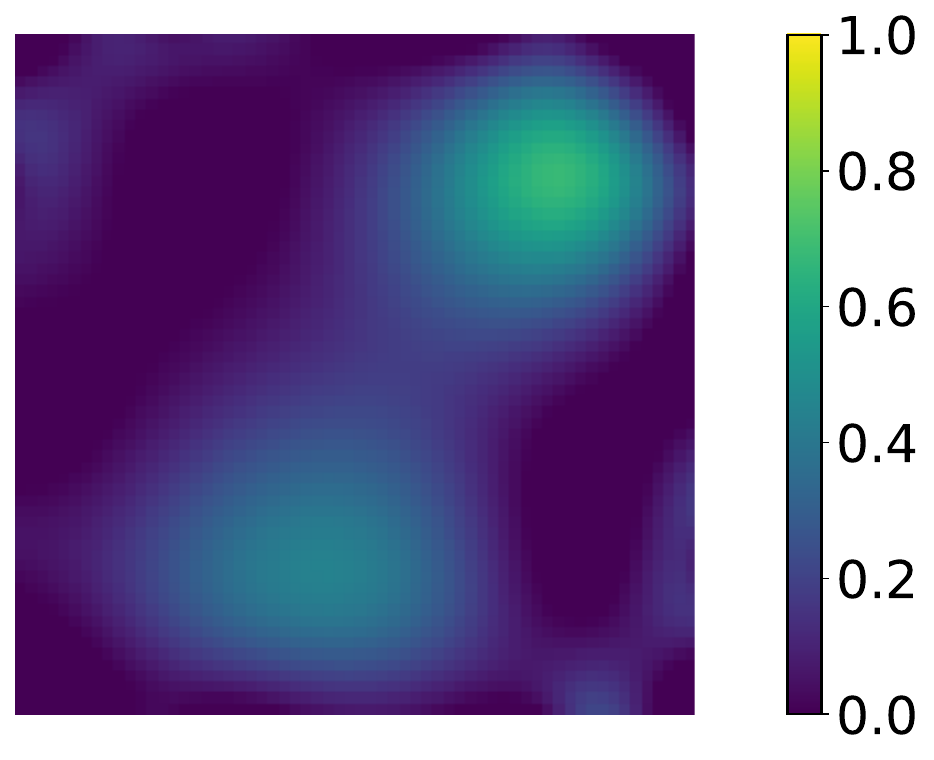} &
        \includegraphics[width=.21\textwidth]{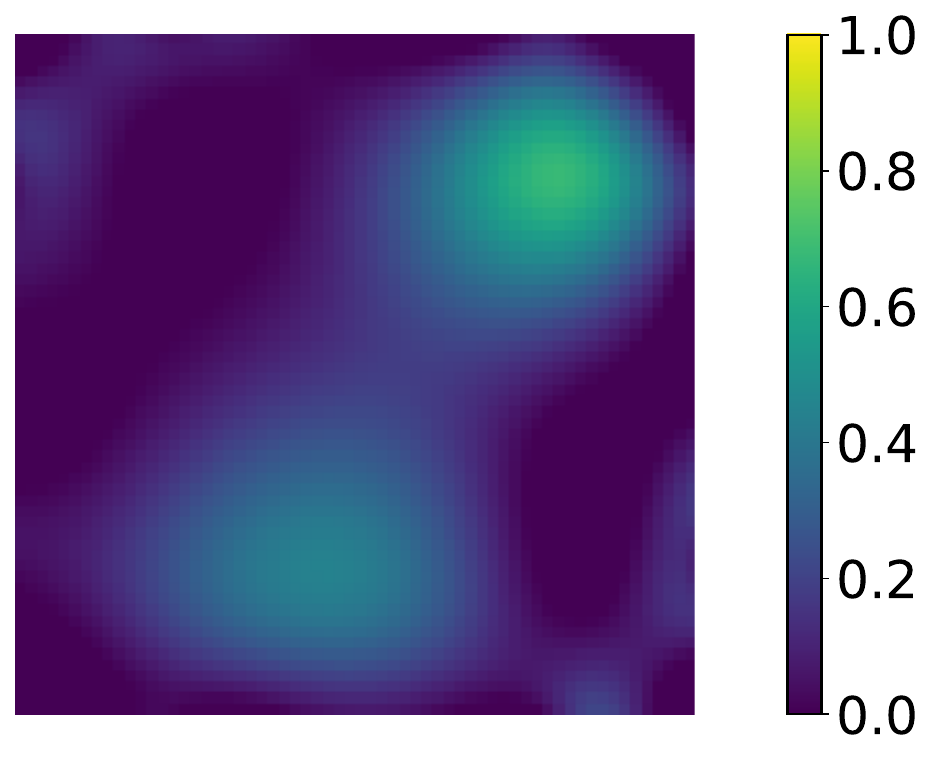} &
        \includegraphics[width=.21\textwidth]{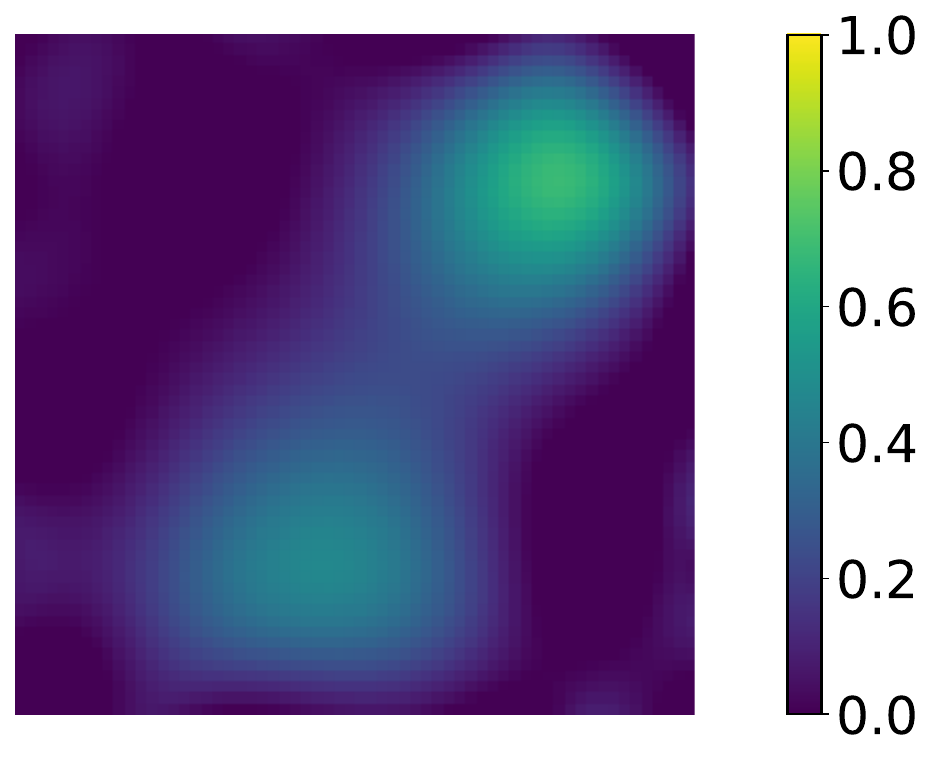} &
        \includegraphics[width=.21\textwidth]{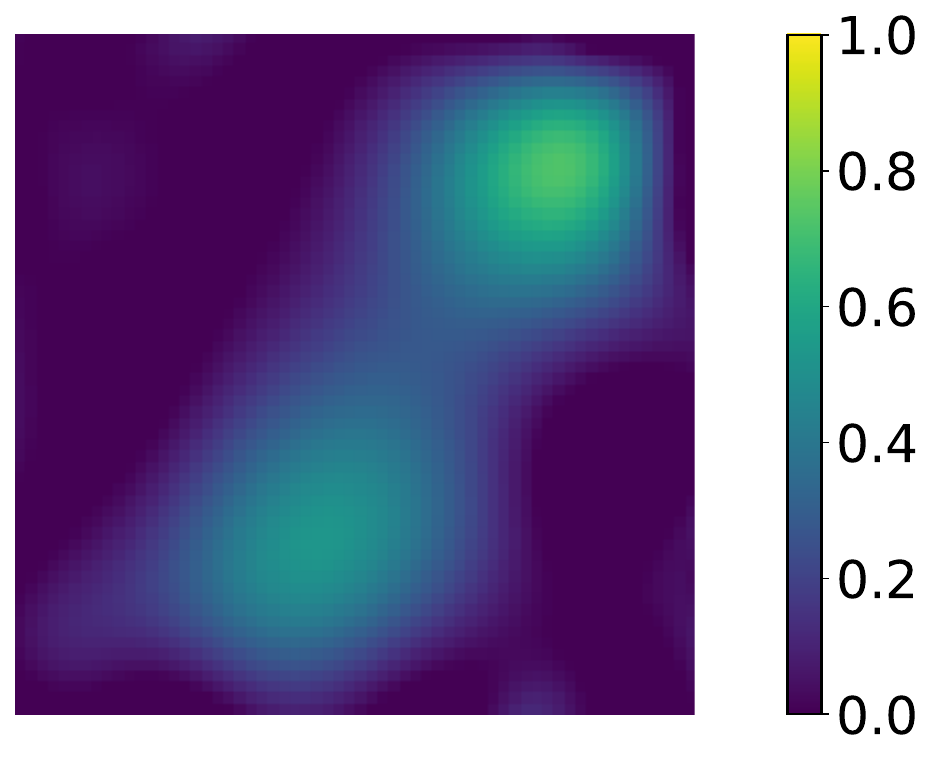} \\
        (a) $1e$-$4$ & (b) $1e$-$8$ & (c) $1e$-$12$ & (d) $1e$-$16$
    \end{tabular}
    \caption{Perturbation reconstructions using different levels of regularization in the MIMO Schrödinger problem.}
    \label{fig:reg:schrodinger}
\end{figure}

Figure~\ref{fig:reg:schrodinger} compares the increased level of regularization for the Schrödinger problem. We observe that artifacts begin to appear for high and low thresholds. However, the quality is consistent for each of the regularization thresholds. 

\begin{figure}[h!]
    \centering
    \begin{tabular}{cccc}
        \includegraphics[width=.21\textwidth]{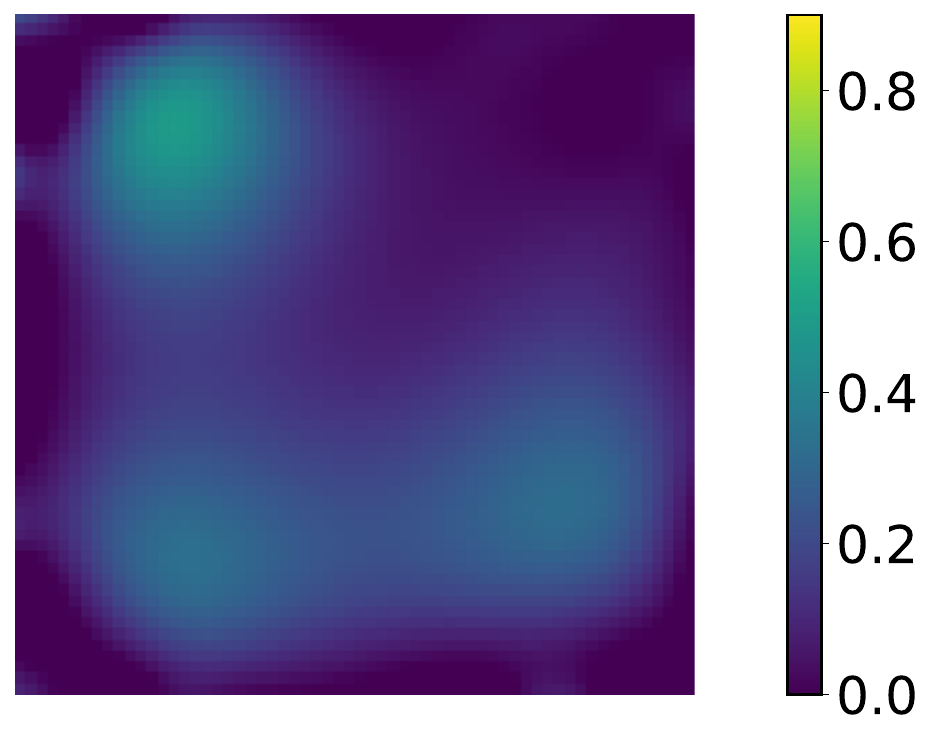} &
        \includegraphics[width=.21\textwidth]{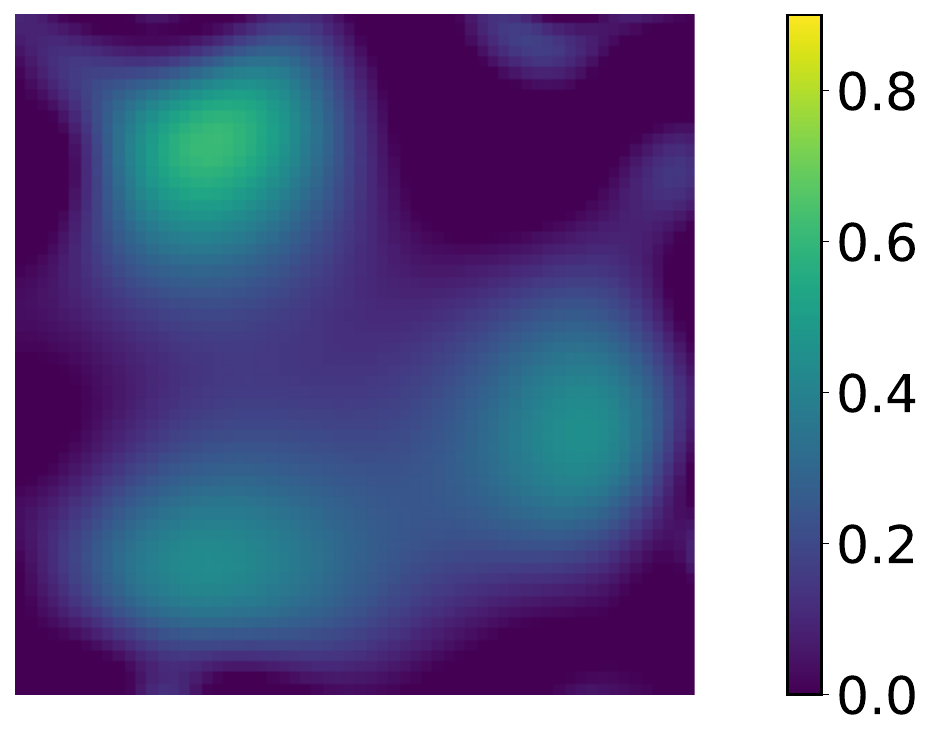} &
        \includegraphics[width=.21\textwidth]{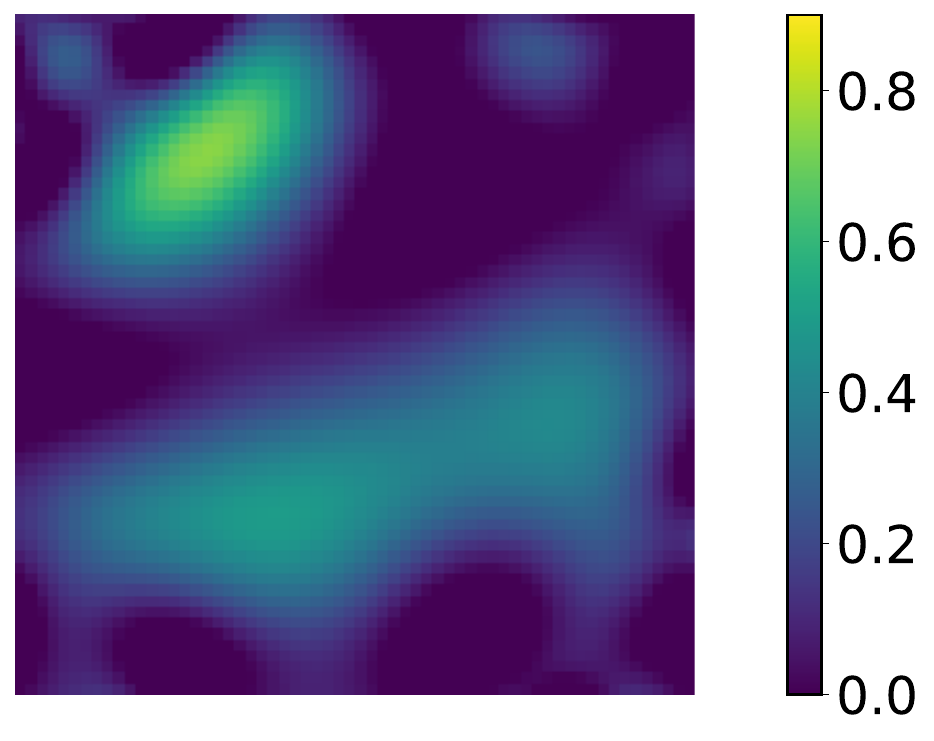} &
        \includegraphics[width=.21\textwidth]{img/helmholtz_reglsl} \\
        (a) $1e$-$4$ & (b) $1e$-$8$ & (c) $1e$-$12$ & (d) $1e$-$16$
    \end{tabular}
    \caption{Perturbation reconstructions using different levels of regularization in the MIMO Helmholtz problem.}
    \label{fig:reg:helmholtz}
\end{figure}

Figure~\ref{fig:reg:helmholtz} compares the increased level of regularization for the Helmholtz problem. For high regularization thresholds, the resolution is quite low. As the threshold decreases, the accuracy of the reconstruction increases, and the resolution improves.

\subsection{Sensitivity Analysis}

In this section, we add noise to the data and analyze the sensitivity of the reconstruction to the noise. The perturbation and selection of frequencies are the same as above. The level of truncation for the mass matrix is fixed to $5e$-$14$ for both Schrödinger and Helmholtz problems. The truncation for the linear inverse problem with $1\%,2\%$ and $5\%$ noise was chosen to be $5e$-$4$, $8e$-$3$ and $4e$-$3$ for the Schrödinger problem and a fixed $5e$-$4$ for the Helmholtz problem.
The noise is sampled from the uniform distribution with a range of $[-1,1]$ and centered at $0$. The percentage of noise is the percent scaling of $|F_0-F|$ added to each component of the data.

\begin{figure}[h!]
    \centering
    \begin{tabular}{cccc}
        \includegraphics[width=.21\textwidth]{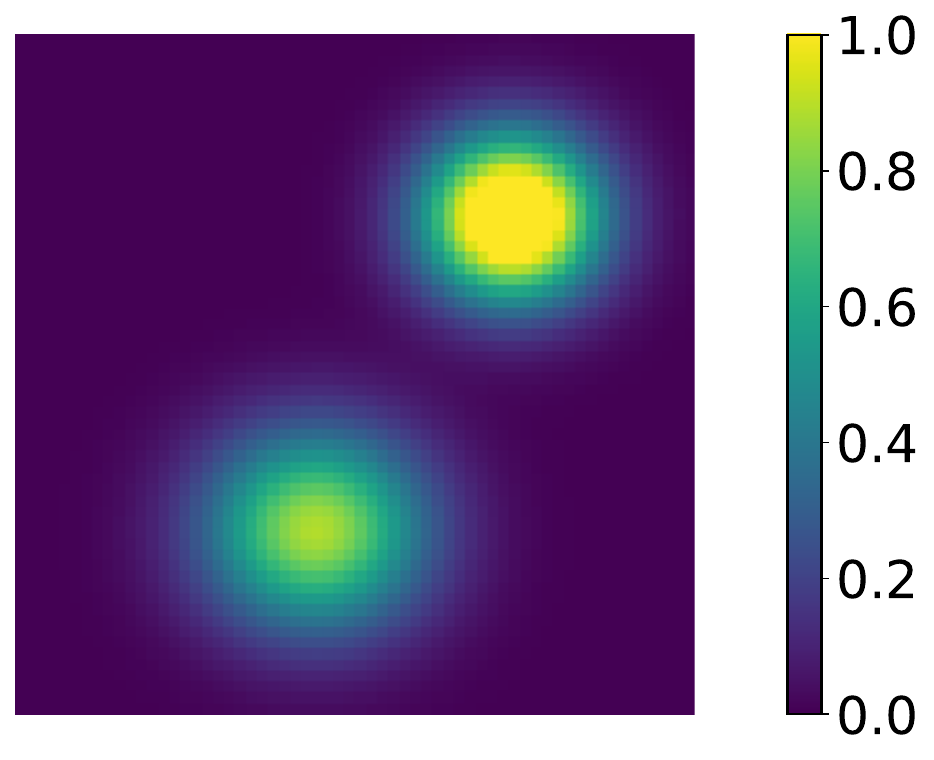} &
        \includegraphics[width=.21\textwidth]{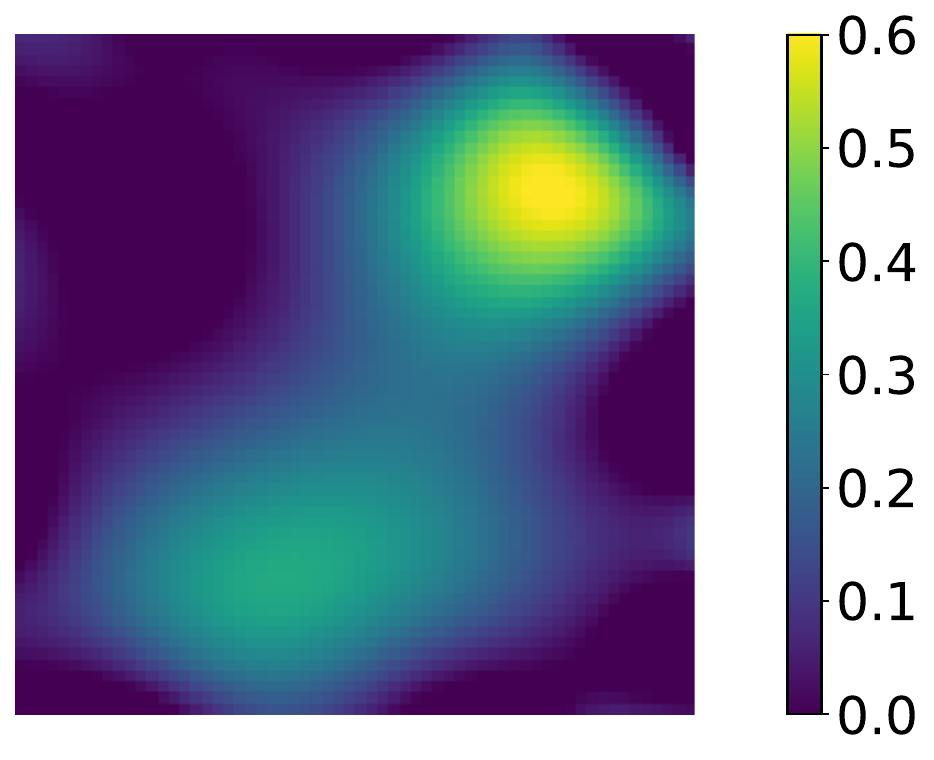} &
        \includegraphics[width=.21\textwidth]{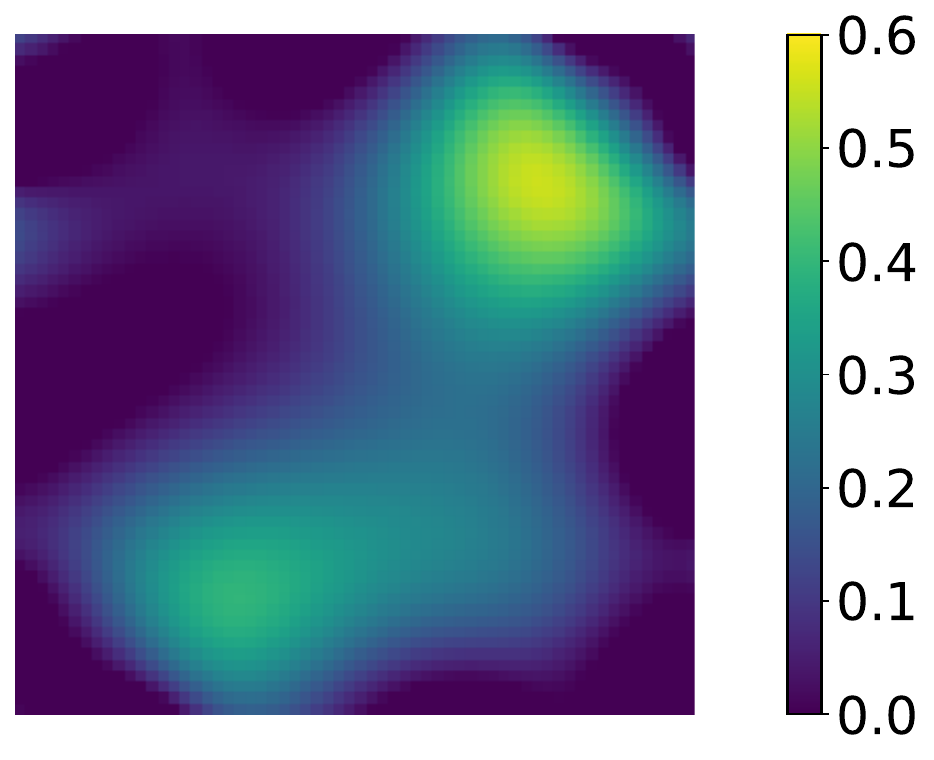} &
        \includegraphics[width=.21\textwidth]{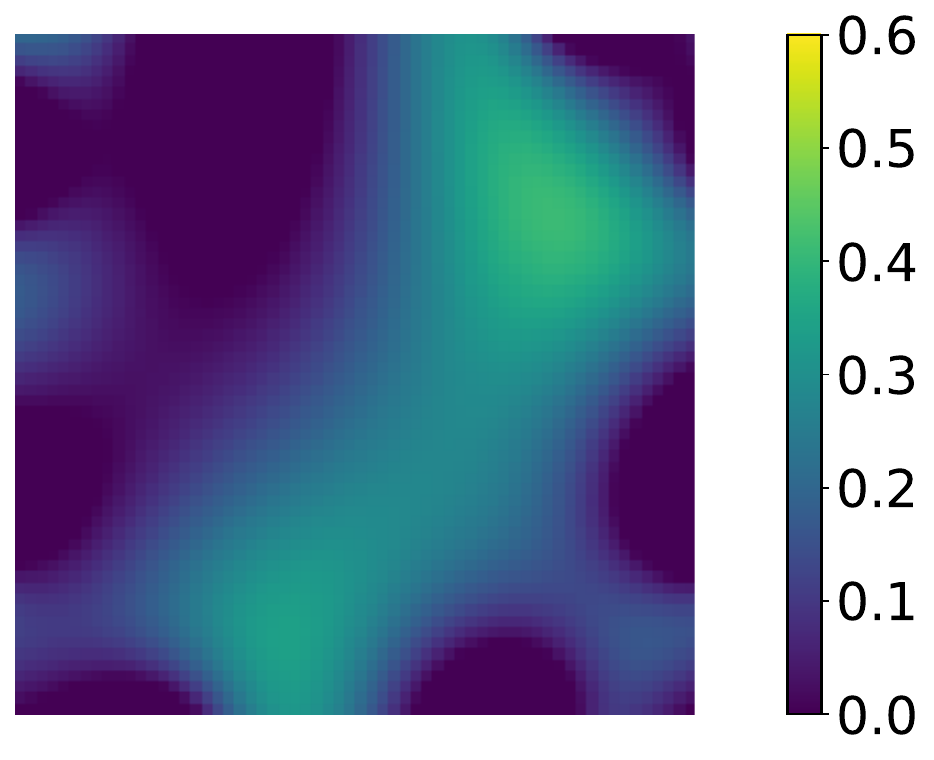} \\
        (a) True & (b) 1\% Noise & (c) 2\% Noise & (d) 5\% Noise
    \end{tabular}
    \caption{Perturbation reconstructions from noisy data in the MIMO Schrödinger problem.}
    \label{fig:noise:schrodinger}
\end{figure}

Figure~\ref{fig:noise:schrodinger} compares the reconstruction of the perturbation as the noise level is increased in the MIMO Schrödinger problem.
We observe that as the level of noise increases, there are significant artifacts in the image of the perturbation.
However, the Reg-LSL algorithm is able to reconstruct the perturbation well until $5\%$ noise is added.

\begin{figure}[h!]
    \centering
    \begin{tabular}{cccc}
        \includegraphics[width=.21\textwidth]{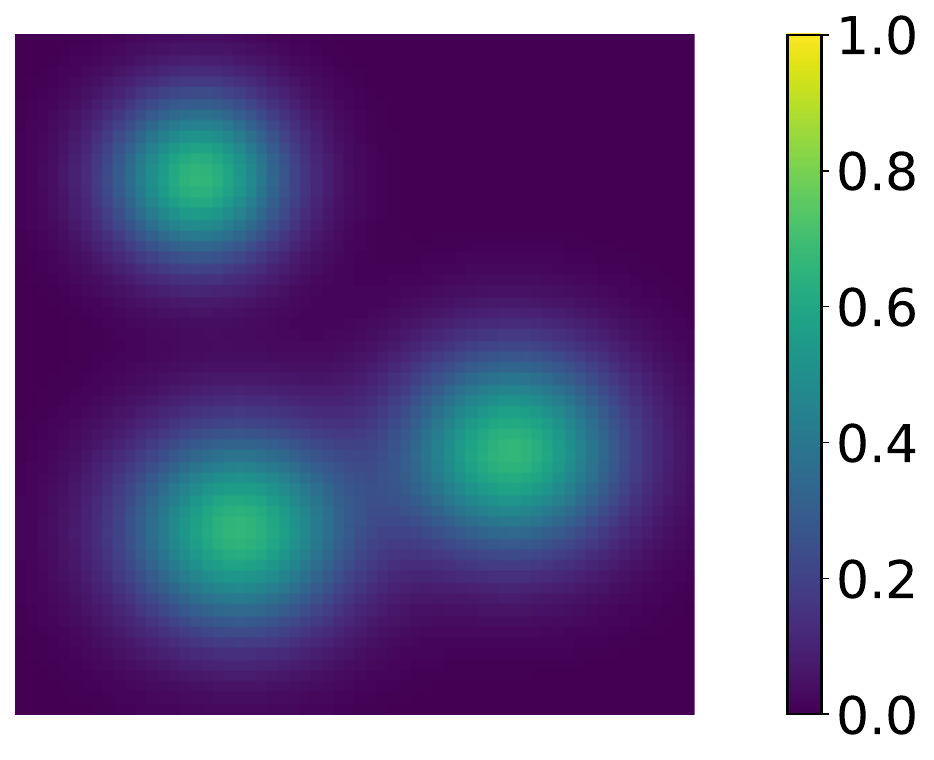} &
        \includegraphics[width=.21\textwidth]{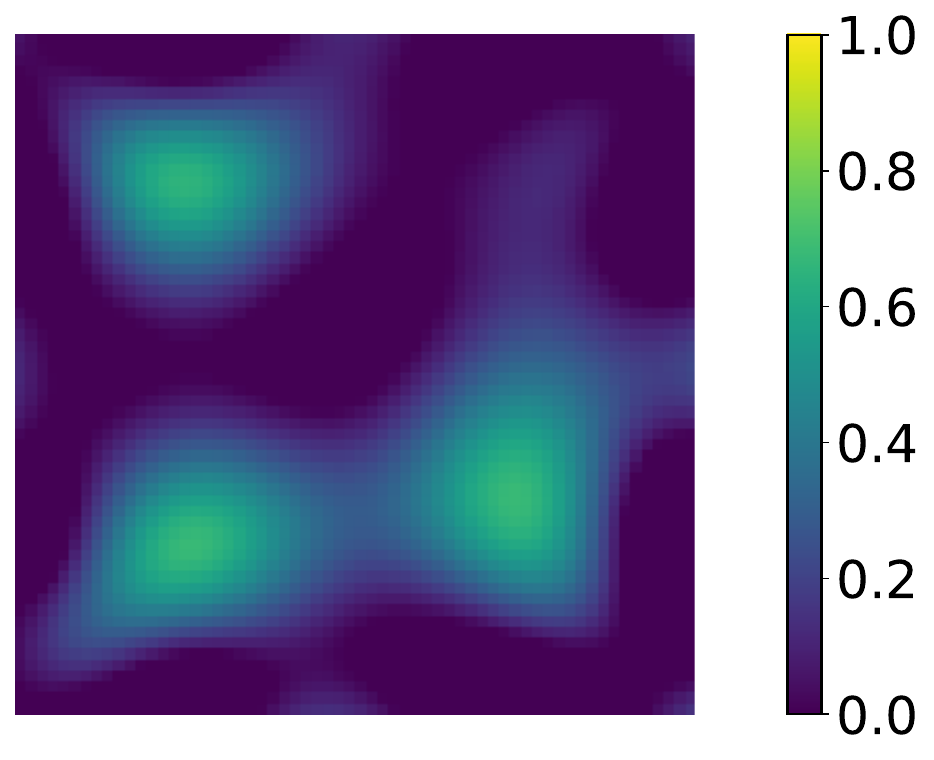} &
        \includegraphics[width=.21\textwidth]{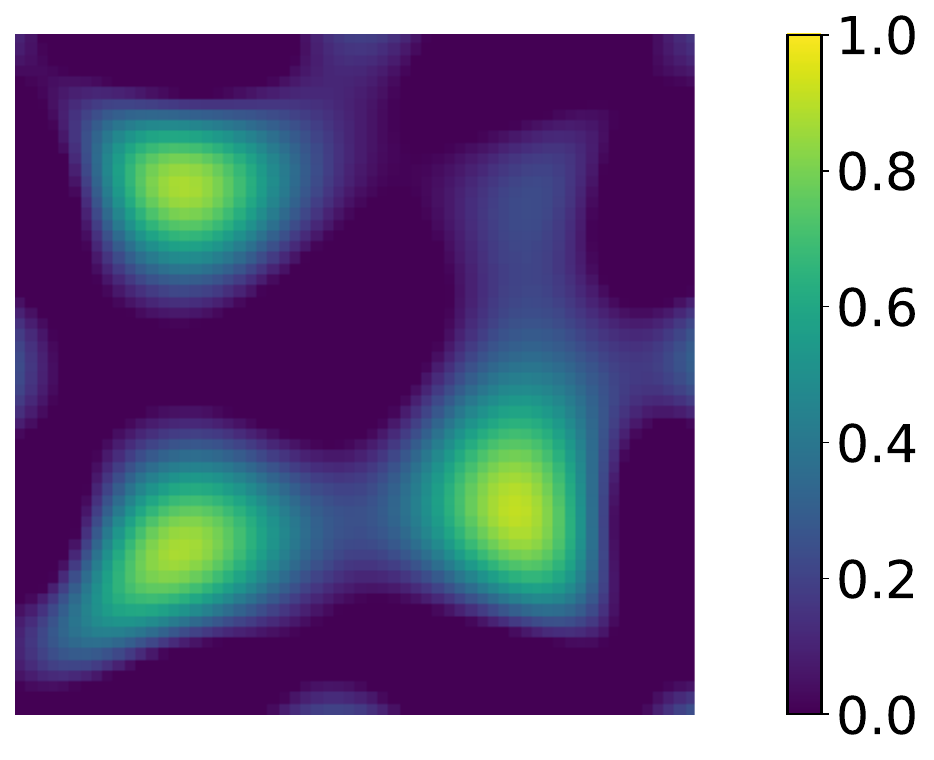} &
        \includegraphics[width=.21\textwidth]{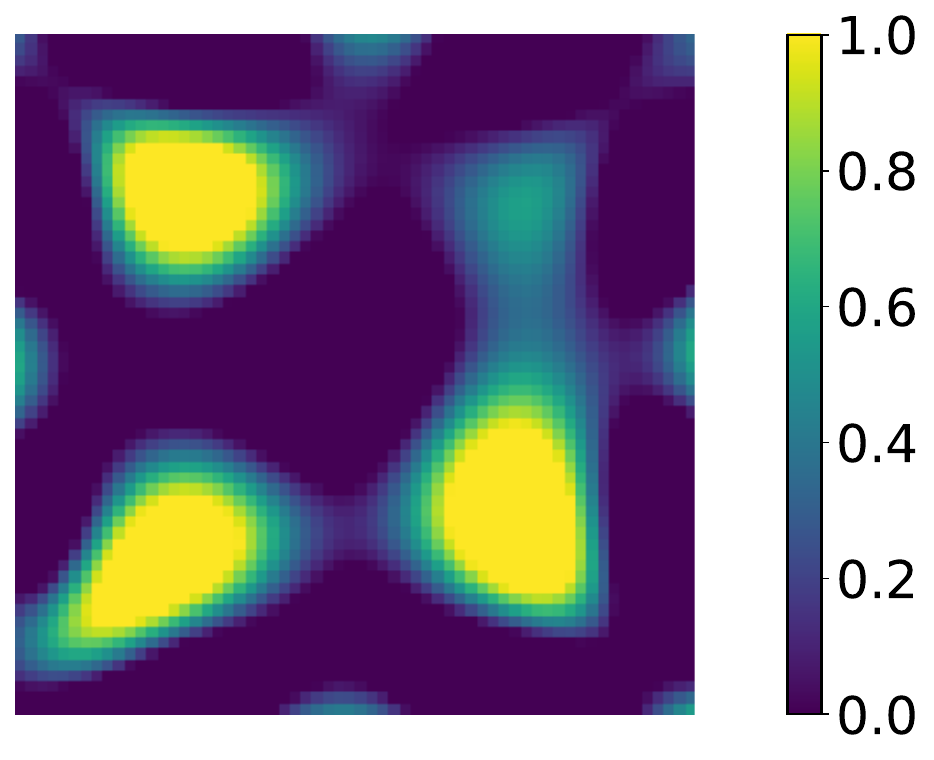} \\
        (a) True & (b) 1\% Noise & (c) 2\% Noise & (d) 5\% Noise
    \end{tabular}
    \caption{Perturbation reconstruction from noisy data in the MIMO Helmholtz problem.}
    \label{fig:noise:helmholtz}
\end{figure}

Figure~\ref{fig:noise:helmholtz} illustrates the reconstructed perturbation using Reg-LSL when noise is added to the data in the MIMO Helmholtz problem.

\section{Conclusion}
In this work, we introduced a double regularization procedure for the Lippmann-Schwinger-Lanczos method in the frequency domain. This includes a data driven truncation of the mass matrix along with the regularization of the linear integral equation. This allows us to improve reconstructions with larger data sets and handle noisy data. We confirmed this with several numerical experiments. 

We also showed here that the Lanczos orthogonalization of the data driven ROM to orthogonalize the frequency domain solutions corresponds precisely to sequential Gram-Schmidt on the corresponding ROM time snapshots. This explains the weak dependence of the orthogonalized snapshots on the medium perturbations for the Schrodinger problem and suggests that the LSL method may work directly for more general formulations. 
 We analyzed several cases of localized conductivity perturbations (corresponding to lower wave speed) in our numerical experiments and showed that while the data-generated internal solutions deviate from the true solutions, the direct reconstructions can still be of high quality.

\section{Acknowledgement}
The authors 
gratefully acknowledge support from the Division of Mathematical Sciences at the US National Science Foundation (NSF) through grants 
DMS-2008441, DMS-2111117, and DMS-2110773,  and from the Air Force Office of Scientific Research
(AFOSR) through grants FA955020-1-0079 and FA9550-23-1-0220. 

~



\bibliography{references.bib}





\end{document}

%% file: math-commands.tex
\usepackage{amsfonts,bm}

\def\eqref#1{equation~(\ref{#1})}

\def\1{\bm{1}}

\def\vb{{\bm{b}}}
\def\vc{{\bm{c}}}
\def\vd{{\bm{d}}}

\def\mF{{\bm{F}}}

\def\mL{{\bm{L}}}
\def\mM{{\bm{M}}}

\def\mQ{{\bm{Q}}}

\def\mS{{\bm{S}}}
\def\mT{{\bm{T}}}

\DeclareMathAlphabet{\mathsfit}{\encodingdefault}{\sfdefault}{m}{sl}
\SetMathAlphabet{\mathsfit}{bold}{\encodingdefault}{\sfdefault}{bx}{n}
\newcommand{\tens}[1]{\bm{\mathsfit{#1}}}

\def\tB{{\tens{B}}}
\def\tC{{\tens{C}}}

\def\tF{{\tens{F}}}

\def\tM{{\tens{M}}}

\def\tQ{{\tens{Q}}}

\def\tS{{\tens{S}}}

